\documentclass[11pt]{amsart}

\usepackage[T1]{fontenc}
\usepackage{lmodern}
\usepackage{amsmath,amssymb,amsthm,mathtools}
\usepackage{microtype}
\usepackage[hidelinks]{hyperref}

\newtheorem{theorem}{Theorem}[section]
\newtheorem{lemma}[theorem]{Lemma}
\newtheorem{proposition}[theorem]{Proposition}
\newtheorem{corollary}[theorem]{Corollary}
\theoremstyle{remark}
\newtheorem{remark}[theorem]{Remark}
\newtheorem{example}[theorem]{Example}
\newtheorem*{counterexample}{Counterexample}

\numberwithin{equation}{section}

\newcommand{\R}{\mathbb{R}}
\newcommand{\Sph}{\mathbb{S}}
\newcommand{\Om}{\Omega}
\newcommand{\dd}{\,\mathrm{d}}
\newcommand{\abs}[1]{\lvert #1\rvert}

\title[Hardy weights from local volume]{Sharp bounds for Hardy weights in terms of local volume}

\author{Tserendorj Batbold}
\address{Department of Mathematics, National University of Mongolia, Ulaanbaatar 14201, Mongolia}
\email{tsbatbold@hotmail.com}
\hypersetup{
  pdftitle={Sharp bounds for Hardy weights in terms of local volume},
  pdfauthor={Tserendorj Batbold}
}


\begin{document}
\begin{abstract}
In dimensions $d\ge2$, we give an explicit two-branch formula for
the least possible Hardy weight at the center of a ball when the
fraction of the ball occupied by an open set is prescribed. If the
volume fractions are prescribed at all radii about the same center,
a single open set symmetric about that center minimizes the weight
there for every power $\alpha>0$ and
every shift $\ell\ge0$. We also derive an integral lower bound that,
on every bounded centrally symmetric convex domain, equals the minimum
of the weight over the domain.

These estimates yield spectral lower bounds for the Robin Laplacian,
Dirichlet polyharmonic operators, and the Dirichlet sub-Laplacian on
the Heisenberg group. We compare the bounds for the Robin and Dirichlet
Laplacians with earlier estimates.
\end{abstract}

\maketitle

\section{Introduction}

Suppose that an open set $\Om\subset\R^d$, $d\ge2$, occupies a
fraction $q$ of a ball $B_r(x)$ centered at a point $x\in\Om$.
The concentric ball of radius $rq^{1/d}$ occupies the same fraction.
Does it minimize the Hardy weight at $x$ among all open sets with
this property? We show that it does up to a threshold in $q$ and
describe a different extremizer above that threshold.

For $\omega\in\Sph^{d-1}$, let
\begin{equation}
\delta_\omega(x):=\inf\{\abs{s}:s\in\R,\ x+s\omega\notin\Om\}
\end{equation}
be the distance from $x$ to the complement along the line
$x+\R\omega$, with $\inf\varnothing=\infty$. Write
\begin{equation}
\psi_r(x):=\frac{\abs{\Om\cap B_r(x)}}{\abs{B_r(x)}}.
\end{equation}
For $\alpha>0$ and $\ell\ge0$, the weight to be minimized is
\begin{equation}
\mathcal W_{\alpha,\ell}(x)
:=\frac{1}{\abs{\Sph^{d-1}}}\int_{\Sph^{d-1}}
   \bigl(\delta_\omega(x)+\ell\bigr)^{-\alpha}\,\dd\omega,
\label{eq:shifted-weight}
\end{equation}
where an infinite directional distance contributes zero.
The case $\alpha=2$, $\ell=0$ occurs in Davies's Hardy
inequality~\cite{Davies1984} and in his local-volume estimate for
the Dirichlet Laplacian~\cite[Theorem~1.5.11]{Davies1989}.
If $\delta_\omega(x)>r$, the entire segment from $x$ to $x+r\omega$
lies in $\Om$. Using this observation, Frank and
Larson~\cite[Lemma~1]{FrankLarson2021} showed that
$r^2\mathcal W_{2,0}(x)\ge 1-\psi_r(x)$.
Estimates for other powers and for positive shifts are given in
\cite[Lemma~3.1]{FrankLaptevSuragan2026}.

At a fixed radius, the sharp bound has two branches. On the first,
the concentric ball gives the minimum $(rq^{1/d}+\ell)^{-\alpha}$.
On the second, an extremizing set contains the whole line through
$x$ in some directions, while in the remaining directions the
complement is at a common finite distance from $x$. The minimum on
this branch is proportional to $1-q$. For $\ell=0$ the two branches
meet at $q=\alpha/(d+\alpha)$.

Now prescribe $\psi_r(x)$ at every radius about the same point $x$.
These data determine the exterior volume fraction of each concentric
annulus. Among all open sets with these same ball fractions, the
least possible proportion of directions with $\delta_\omega(x)<t$
equals the supremum of the exterior fractions of annuli with outer
radius at most $t$. One open set symmetric about $x$ attains this
minimum for every $t>0$. Since each Hardy weight is an integral of
these proportions against a positive kernel, the same set minimizes
$\mathcal W_{\alpha,\ell}(x)$ for all $\alpha>0$ and $\ell\ge0$.

For spectral applications we need estimates that hold uniformly
in $x$. Put $\Psi_r:=\sup_{x\in\Om}\psi_r(x)$. The full function
$r\mapsto\Psi_r$ gives an integral lower bound for
$\inf_{x\in\Om}\mathcal W_{\alpha,\ell}(x)$. On bounded centrally
symmetric convex domains this bound is exact. On every such domain
other than a ball, it strictly improves the bound obtained by
optimizing the fixed-radius estimate over $r$.

Section~\ref{sec:main-results} states the main geometric results.
Section~\ref{sec:unshifted-proof} proves the fixed-radius estimates.
Section~\ref{sec:shifted} treats the full volume profile and derives
the integral bound.

For the first Dirichlet eigenvalue, Lieb's intersection
theorem~\cite{Lieb1983} complements the Faber--Krahn bound in terms
of total volume with a lower bound in terms of ball intersections.
Lieb allows centers in all of $\R^d$. With our choice of centers
the corresponding generalized inradius at level $q\in(0,1)$ is
\[
R_\Om^{(q)}:=\sup\{r>0:\Psi_r\ge q\}.
\]
A generalized inradius defined by volume fractions also enters
semiclassical estimates for eigenvalue Riesz
means~\cite{FrankLarsonPfeiffer2026}. Measuring the part of a ball
outside $\Om$ by capacity leads to the interior capacitary radius
and to two-sided spectral bounds~\cite{MazyaShubin2005}.

For the Robin Laplacian, the boundary parameter enters through the
shift in the weight. Suppose that $\Om$ has uniformly Lipschitz
boundary and that $\sigma>0$ is constant. Write
$\lambda_\Om^{\mathrm{R}\sigma}$ for the bottom of the spectrum of
the operator associated with the quadratic form
\begin{equation}
h_\Om^{\mathrm{R}\sigma}[u]
:=\int_\Om |\nabla u|^2\,\dd x
  +\sigma\int_{\partial\Om}|u|^2\,\dd S,
\qquad u\in H^1(\Om).
\label{eq:robin-form}
\end{equation}
Here the boundary integral uses the Sobolev trace of $u$ and surface
measure $\dd S$. On bounded domains satisfying their strong local
Lipschitz condition, Kova\v{r}\'{i}k and
Laptev~\cite[Theorem~4.2]{KovarikLaptev2012} proved that
\[
h_\Om^{\mathrm{R}\sigma}[u]\ge
\frac14\int_\Om \mathcal W_{2,(2\sigma)^{-1}}(x)|u(x)|^2\,\dd x.
\]
For open sets with uniformly Lipschitz boundary, Theorem~2.1 in
the \href{https://arxiv.org/abs/2509.18878v1}{arXiv version} of
Frank, Laptev, and Suragan~\cite{FrankLaptevSuragan2026} states that
\begin{equation}
\lambda_\Om^{\mathrm{R}\sigma}
\ge \frac{d\sigma^2}{(1+2\sigma r)^2}(1-\Psi_r),
\qquad r>0.
\label{eq:earlier-robin-bound}
\end{equation}
Unfortunately, \eqref{eq:earlier-robin-bound} can fail even on a
bounded rectangular box.

\begin{counterexample}
Let $d=2049$, $r=1$, $\sigma=1/2$, and
\[
\Om=(-1/96,1/96)\times(-64,64)^{2048}.
\]
Then $\lambda_\Om^{\mathrm{R}\sigma}\le64$, whereas the right side
of \eqref{eq:earlier-robin-bound} is at least $2049/32$.
\end{counterexample}

\begin{proof}
The constant function is an admissible trial function on this box,
so $\lambda_\Om^{\mathrm{R}\sigma}\le\sigma|\partial\Om|/|\Om|$.
For a box with side lengths $2a_j$, the two faces perpendicular to
the $j$th coordinate axis have combined area $|\Om|/a_j$. Thus
\[
\lambda_\Om^{\mathrm{R}\sigma}
\le \frac12\left(96+\frac{2048}{64}\right)=64.
\]

To estimate $\Psi_1$, denote by $v_k$ the volume of the unit ball
in $\R^k$. Every section of $B_1(x)$ perpendicular to the first
coordinate axis has $(d-1)$-dimensional volume at most $v_{d-1}$.
Since $\Om$ has width $1/48$ in this direction, Fubini's theorem
gives $\abs{\Om\cap B_1(x)}\le v_{d-1}/48$ for every center $x$.
Slicing the unit ball in the same direction gives
\[
\frac{v_d}{v_{d-1}}
=\int_{-1}^{1}(1-t^2)^{1024}\,\dd t
\ge \int_{-1/32}^{1/32}(1-1024t^2)\,\dd t
=\frac1{24}.
\]
Here $(d-1)/2=1024$, and the inequality follows from Bernoulli's
inequality on the interval where $1-1024t^2\ge0$. Consequently
$\Psi_1\le v_{d-1}/(48v_d)\le1/2$, and
\[
\frac{d\sigma^2}{(1+2\sigma r)^2}(1-\Psi_r)
=\frac{2049}{16}(1-\Psi_1)
\ge\frac{2049}{32}>64.
\]
\end{proof}

We return to the Robin problem in Section~\ref{sec:applications},
where we combine a Hardy inequality for \eqref{eq:robin-form} with
our geometric estimates to obtain lower bounds for
$\lambda_\Om^{\mathrm{R}\sigma}$ in terms of local volume.
The same section gives applications to Dirichlet polyharmonic
operators and to the Dirichlet sub-Laplacian on the Heisenberg group.
Section~\ref{sec:comparisons} compares the resulting spectral bounds
with earlier estimates.

\section{Main results}
\label{sec:main-results}

\subsection{Sharp bounds at a prescribed radius}

Theorem~\ref{thm:shifted-geometric} gives the sharp lower bound for
every $\ell\ge0$ in dimensions $d\ge2$. We first state the explicit
formula at zero shift. The corresponding lower bound also holds in
dimension one.

\begin{theorem}
\label{thm:geometric}
Let $d\in\mathbb N$, let $\Om\subset\R^d$ be nonempty and open, and
let $\alpha>0$. For every $x\in\Om$ and $r>0$,
\begin{equation}
\frac{1}{\abs{\Sph^{d-1}}}\int_{\Sph^{d-1}}
\delta_\omega(x)^{-\alpha}\,\dd\omega
\ge \frac{1}{r^\alpha}\Phi_{\alpha,d}\bigl(\psi_r(x)\bigr).
\label{eq:main-geom-intro}
\end{equation}
Set
\begin{equation}
q_*:=\frac{\alpha}{d+\alpha},
\qquad
K_{\alpha,d}:=\frac{q_*^{-\alpha/d}}{1-q_*}
=\frac{d+\alpha}{d}\left(\frac{d+\alpha}{\alpha}\right)^{\alpha/d}.
\label{eq:K-alpha-d}
\end{equation}
Define
\begin{equation}
\Phi_{\alpha,d}(q)=
\begin{cases}
q^{-\alpha/d}, & 0<q\le q_*,\\[5pt]
K_{\alpha,d}(1-q), & q_*\le q\le1.
\end{cases}
\label{eq:Phi}
\end{equation}
\end{theorem}

For $d\ge2$, Theorem~\ref{thm:sharpness} gives open sets realizing
equality for every $q\in(0,1]$ and shows that bounded open sets give
the same infimum.

For the shifted bound, let $d\ge2$, $\alpha>0$, and $\ell\ge0$.
We keep $\ell$ fixed when varying $r$. For each $r>0$, let $T_r$
be the unique solution in $(0,r)$ of
\begin{equation}
(d+\alpha)T_r^d+d\ell T_r^{d-1}=\alpha r^d.
\label{eq:shifted-tangency}
\end{equation}
The polynomial on the left is continuous and strictly increasing on
$[0,\infty)$. It starts at zero and exceeds $\alpha r^d$ at $r$, which
gives existence and uniqueness. Define the threshold volume fraction
and the linear coefficient by
\begin{equation}
q_r^*:=\left(\frac{T_r}{r}\right)^d,
\qquad
C_{\alpha,d,\ell}(r):=\frac{(T_r+\ell)^{-\alpha}}{1-q_r^*}.
\label{eq:shifted-linear-coef}
\end{equation}
The sharp lower bound at radius $r$ is given by
\begin{equation}
G_{\alpha,d,\ell,r}(q):=
\begin{cases}
(r q^{1/d}+\ell)^{-\alpha}, & 0<q\le q_r^*,\\[5pt]
C_{\alpha,d,\ell}(r)(1-q), & q_r^*\le q\le1.
\end{cases}
\label{eq:shifted-profile}
\end{equation}
The two branches agree at $q_r^*$, and each decreases with $q$.
Thus $G_{\alpha,d,\ell,r}$ is continuous and decreasing on $(0,1]$.
At $\ell=0$, equation
\eqref{eq:shifted-tangency} gives $T_r/r=(\alpha/(d+\alpha))^{1/d}$.
With $K_{\alpha,d}$ as in \eqref{eq:K-alpha-d}, this gives
\[
C_{\alpha,d,0}(r)=r^{-\alpha}K_{\alpha,d},
\qquad
G_{\alpha,d,0,r}(q)=r^{-\alpha}\Phi_{\alpha,d}(q).
\]

\begin{theorem}
\label{thm:shifted-geometric}
Let $d\ge2$, $\alpha>0$, and $\ell\ge0$. For every nonempty open set
$\Om\subset\R^d$, $x\in\Om$, and $r>0$,
\begin{equation}
\mathcal W_{\alpha,\ell}(x)
\ge G_{\alpha,d,\ell,r}(\psi_r(x)).
\label{eq:shifted-geom}
\end{equation}
For every $r>0$ and $q\in(0,1]$, there is an open set $\Om$ and a
point $x\in\Om$ with $\psi_r(x)=q$ for which equality holds in
\eqref{eq:shifted-geom}. The value $G_{\alpha,d,\ell,r}(q)$ is also
the infimum of $\mathcal W_{\alpha,\ell}(x)$ over bounded open sets
and points with $\psi_r(x)=q$.
\end{theorem}

\subsection{Volume fractions at all radii}

Fix $x\in\Om$. Knowing $\psi_r(x)$ for all $r>0$ determines the
fraction of each concentric annulus outside the set. Write
$D_x(r)=r^d(1-\psi_r(x))$ for $r>0$ and set $D_x(0)=0$.
For $t>0$, define
\begin{equation}
H_x(t)=\sup_{0\le a<b\le t}
\frac{D_x(b)-D_x(a)}{b^d-a^d}.
\label{eq:annular-deficit}
\end{equation}
The quotient is the fraction of $B_b(x)\setminus B_a(x)$ outside
$\Om$, with $B_0(x)=\varnothing$.

\begin{theorem}
\label{thm:all-radii-sharp}
Let $d\ge2$, let $\Om\subset\R^d$ be open, and fix $x\in\Om$.
Let $\mu$ denote normalized surface measure on $\Sph^{d-1}$.
For an open set $\Om'$ containing $x$, write $\psi_s^{\Om'}(x)$
and $\delta_\omega^{\Om'}(x)$ for its volume fractions and directional
distances. If $\psi_s^{\Om'}(x)=\psi_s^{\Om}(x)$ for every $s>0$, then
\begin{equation}
H_x(t)\le\mu\{\omega:\delta_\omega^{\Om'}(x)<t\},\qquad t>0.
\label{eq:annular-tail}
\end{equation}
There is an open set $\Om'$ symmetric about $x$, with the same
volume fractions, for which equality holds for every $t>0$. It can
be taken bounded if $\Om$ is bounded.

Consequently, among those open sets the minimum of the Hardy weight
at $x$ is
\begin{equation}
\mathcal L_{\alpha,\ell}(\Om,x)
:=\alpha\int_0^\infty H_x(t)(t+\ell)^{-\alpha-1}\,\dd t,
\qquad \alpha>0,\quad \ell\ge0.
\label{eq:all-radii-sharp}
\end{equation}
The same set attains the minimum simultaneously for all $\alpha$ and
$\ell$.
\end{theorem}

Fix $r>0$ and $q\in(0,1]$. For each realizable profile with
$\psi_r(x)=q$, the sets that share that profile satisfy the constraint
in Theorem~\ref{thm:shifted-geometric}. Their minimum in
\eqref{eq:all-radii-sharp} is therefore at least $G_{\alpha,d,\ell,r}(q)$.
Take now the profile of a set attaining the fixed-radius bound.
Since this set has weight $G_{\alpha,d,\ell,r}(q)$, the minimum for its
profile is at most this value. Thus minimizing over all profiles
with $\psi_r(x)=q$ recovers the fixed-radius minimum.

When only the suprema $\Psi_r$ are known, Proposition~\ref{prop:all-radii}
gives an explicit integral lower bound for
$\inf_{x\in\Om}\mathcal W_{\alpha,\ell}(x)$.
Proposition~\ref{prop:symmetric-convex} shows that on bounded centrally
symmetric convex sets this integral gives the exact minimum over
centers and, except on balls, strictly exceeds the optimized
fixed-radius bound. This integral gives spectral lower bounds in
Section~\ref{sec:applications}.

\section{Sharp bounds at a fixed radius}
\label{sec:unshifted-proof}

Fix $x\in\Om$ and $r>0$, and put $Y(\omega)=\delta_\omega(x)/r$.
A ball centered at $x$ lies in $\Om$, so its radius gives a positive
lower bound for $\delta_\omega(x)$ in every direction. For $t>0$,
the inequality $\delta_\omega(x)>t$ is equivalent to
$x+[-t,t]\omega\subset\Om$. The compact segment stays inside $\Om$
for nearby directions. Hence $\omega\mapsto\delta_\omega(x)$ is lower
semicontinuous and $Y$ is measurable.

To relate these distances to the volume fraction, note that
\[
0\le\rho<\delta_\omega(x)
\quad\Longrightarrow\quad
x+\rho\omega\in\Om.
\]
In polar coordinates, the radial integral for $\abs{\Om\cap B_r(x)}$
therefore contains the interval from $0$ to
$\min\{r,\delta_\omega(x)\}$. Thus
\begin{align}
\psi_r(x)
&=
\frac{
\int_{\Sph^{d-1}}\int_0^r
\mathbf 1_{\Om}(x+\rho\omega)\rho^{d-1}\,\dd\rho\,\dd\omega
}{
\int_{\Sph^{d-1}}\int_0^r
\rho^{d-1}\,\dd\rho\,\dd\omega
}
\notag\\
&\ge
\frac{
\int_{\Sph^{d-1}}\int_0^{\min\{r,\delta_\omega(x)\}}
\rho^{d-1}\,\dd\rho\,\dd\omega
}{
\int_{\Sph^{d-1}}\int_0^r
\rho^{d-1}\,\dd\rho\,\dd\omega
}
\notag\\
&=
\frac{1}{|\Sph^{d-1}|}\int_{\Sph^{d-1}}
\min\left\{
1,
\left(\frac{\delta_\omega(x)}{r}\right)^d
\right\}\,\dd\omega.
\label{eq:polar-occupation}
\end{align}
Equip $\Sph^{d-1}$ with the probability measure
\[
\dd\mu(\omega)=\frac{\dd\omega}{|\Sph^{d-1}|}
\]
so that
$\int_{\Sph^{d-1}}\min\{1,Y^d\}\,\dd\mu\le\psi_r(x)$.
Since $r^\alpha\mathcal W_{\alpha,0}(x)$ is the mean of $Y^{-\alpha}$,
we first determine the least possible value of this mean under the
moment constraint above.

\subsection{A moment inequality}

\begin{lemma}
\label{lem:convex-envelope}
Let $(X,\mu)$ be a probability space, let $d,\alpha>0$, and let
$Y:X\to(0,\infty]$ be measurable. We use the conventions
$\infty^{-\alpha}=0$ and $\min\{1,\infty\}=1$.
Suppose that $q\in(0,1]$ satisfies
\begin{equation}
\int_X \min\{1,Y^d\}\,\dd\mu \le q.
\label{eq:moment-constraint}
\end{equation}
Then
\begin{equation}
\int_X Y^{-\alpha}\,\dd\mu
\ge
\Phi_{\alpha,d}(q),
\label{eq:moment-bound}
\end{equation}
where $\Phi_{\alpha,d}$ is given by \eqref{eq:Phi}. For every
$q\in(0,1]$, this value is attained for some probability space and
some $Y$ satisfying \eqref{eq:moment-constraint}.
\end{lemma}

\begin{proof}
To turn \eqref{eq:moment-constraint} into a lower bound on a mean, set
$
p:=\frac{\alpha}{d},
Z:=1-\min\{1,Y^d\}.
$
Since $Y>0$, we have $0\le Z<1$. By definition,
\[
\min\{1,Y^d\}=1-Z.
\]
Since $\mu(X)=1$, integrating this identity in
\eqref{eq:moment-constraint} gives
\begin{equation}
\int_X Z\,\dd\mu\ge 1-q.
\label{eq:EZ}
\end{equation}

For $z\in[0,1)$, let $f(z)$ be the smallest value of $Y^{-\alpha}$
compatible with $Z=z$. This gives $f:[0,1)\to[0,\infty)$ defined by
\begin{equation}
f(0):=0,
\qquad
f(z):=(1-z)^{-p}\quad (0<z<1).
\end{equation}
By construction,
\begin{equation}
Y^{-\alpha}\ge f(Z).
\label{eq:point-f}
\end{equation}
Indeed, if $Y<1$, then $Z>0$ determines $Y$ through $Z=1-Y^d$, and
\[
f(Z)=(Y^d)^{-p}=Y^{-\alpha}.
\]
The value $Z=0$ corresponds to $Y\ge1$, and the choice $Y=\infty$
attains the value zero. In particular, $f(0)=0$ while
$\lim_{z\downarrow0}(1-z)^{-p}=1$. This drop at the endpoint
precludes convexity on $[0,1)$.

To apply Jensen's inequality, we find the greatest convex function
below $f$. A line through $(0,0)$ lies below $f$ if
its slope is at most $f(z)/z$ for all $0<z<1$. Thus its largest
permitted slope is the minimum of this ratio. Set
\begin{equation}
z_*:=\frac{d}{d+\alpha}=\frac{1}{1+p}.
\label{eq:zstar}
\end{equation}
For $0<z<1$,
\[
\frac{\dd}{\dd z}\log\frac{(1-z)^{-p}}{z}
=
\frac{p}{1-z}-\frac1z
=
\frac{(1+p)z-1}{z(1-z)}.
\]
The derivative changes sign from negative to positive at $z_*$,
which is therefore the unique minimizer of the ratio.
The line through $(0,0)$ and $(z_*,f(z_*))$ lies below $f$. The stationarity
condition $z_*f'(z_*)=f(z_*)$ makes it tangent to $f$ at $z_*$. Since
$1-z_*=q_*$, its slope is
\begin{equation}
\frac{(1-z_*)^{-p}}{z_*}
=\frac{q_*^{-\alpha/d}}{1-q_*}=K_{\alpha,d}.
\end{equation}
Because the slopes match at $z_*$ and $f$ is convex on $(0,1)$,
joining the line on $[0,z_*]$ to $f$ on $[z_*,1)$ gives the convex
minorant
\begin{equation}
\varphi(z)
=
\begin{cases}
K_{\alpha,d}z,
&
0\le z\le z_*,
\\[6pt]
(1-z)^{-\alpha/d},
&
z_*\le z<1.
\end{cases}
\label{eq:conv-envelope}
\end{equation}
To prove maximality, let $g$ be any convex minorant of $f$.
For $z\in[0,z_*]$, bound $g(z)$ by the chord joining
$(0,g(0))$ and $(z_*,g(z_*))$. The endpoint inequalities
$g(0)\le f(0)=0$ and $g(z_*)\le f(z_*)$ give
\[
g(z)\le \left(1-\frac{z}{z_*}\right)g(0)+\frac{z}{z_*}g(z_*)
\le \frac{z}{z_*}f(z_*)=\varphi(z).
\]
On $[z_*,1)$, $g\le f=\varphi$, so $\varphi$ is the lower convex
envelope.

By \eqref{eq:point-f} and the definition of $\varphi$,
\[
Y^{-\alpha}\ge f(Z)\ge\varphi(Z).
\]
Since $Z<1$ almost everywhere, its mean belongs to $[0,1)$, the
domain of $\varphi$. Jensen's inequality, the monotonicity of
$\varphi$, and \eqref{eq:EZ} give
\begin{equation}
\int_XY^{-\alpha}\,\dd\mu
\ge
\int_X\varphi(Z)\,\dd\mu
\ge
\varphi\left(\int_XZ\,\dd\mu\right)
\ge
\varphi(1-q).
\end{equation}
If $q\le q_*$, then $1-q\ge z_*$ and we use the curved branch,
\[
\varphi(1-q)=q^{-\alpha/d}.
\]
If $q\ge q_*$, then $1-q\le z_*$ and we use the linear branch,
\[
\varphi(1-q)=K_{\alpha,d}(1-q).
\]
These are the two branches of \eqref{eq:Phi}, so
\eqref{eq:moment-bound} follows.

To attain the bound for $q\le q_*$, take
$Y\equiv q^{1/d}$. The moment constraint is an equality and
$Z=1-q$ is constant. Thus Jensen's
inequality is an equality, and the pointwise inequalities are
equalities because $f(Z)=\varphi(Z)=Y^{-\alpha}$.

For $q>q_*$, let $Z$ take only the values $0$ and $z_*$.
To give $Z$ mean $1-q$, its mass at $z_*$ must be
\[
\theta:=\frac{1-q}{z_*}\in[0,1).
\]
On a two-point probability space take $Y=(1-z_*)^{1/d}$ with
probability $\theta$, so $Z=z_*$, and $Y=\infty$ with probability
$1-\theta$, so $Z=Y^{-\alpha}=0$. At $q=1$ this gives $Y=\infty$
almost surely. Then
\[
\int Z\,\dd\mu=1-q
\]
and the mean of $Y^{-\alpha}$ is
\[
\int Y^{-\alpha}\,\dd\mu
=
\theta(1-z_*)^{-\alpha/d}
=\frac{1-q}{1-q_*}q_*^{-\alpha/d}
=K_{\alpha,d}(1-q).
\]
Thus equality holds in \eqref{eq:moment-constraint}, and the mean of
$Y^{-\alpha}$ attains the linear branch of $\Phi_{\alpha,d}(q)$.
\end{proof}

\begin{remark}
The line is tangent to $f$ at $z=z_*$. Since the bound evaluates
$\varphi$ at $z=1-q$, this contact corresponds to the transition
$q_*=1-z_*=\alpha/(d+\alpha)$ in the volume fraction.
\end{remark}

\subsection{Proof of Theorem~\ref{thm:geometric}}

\begin{proof}
Fix $x\in\Om$ and $r>0$. By \eqref{eq:polar-occupation}, the function
$Y(\omega)=\delta_\omega(x)/r$ satisfies \eqref{eq:moment-constraint}
with respect to normalized surface measure.
Applying Lemma~\ref{lem:convex-envelope} with $q=\psi_r(x)$ yields
\[
\frac{1}{|\Sph^{d-1}|}\int_{\Sph^{d-1}}
\left(\frac{r}{\delta_\omega(x)}\right)^\alpha
\dd\omega
\ge
\Phi_{\alpha,d}(\psi_r(x)).
\]
Dividing by $r^\alpha$ gives \eqref{eq:main-geom-intro}.
\end{proof}

\subsection{Sharpness}

The extremizers in Lemma~\ref{lem:convex-envelope} prescribe a constant
distance or two distances on sets of specified spherical measure. In
dimensions $d\ge2$, open sets realize both configurations.

Fix a radius $r>0$, a distance $a\in(0,r)$, and a fraction
$\theta\in[0,1]$. We construct an open set for which the distance
from $0$ along a fraction $\theta$ of the directions is $a$, with
infinite distance along the remaining directions. Since the
directional distance is invariant under $\omega\mapsto-\omega$,
choose a closed antipodally symmetric set $A\subset\Sph^{d-1}$ of
normalized surface measure $\theta$. Such a set can be chosen from
the family
\[
A_c=\{\omega\in\Sph^{d-1}:|\omega_1|\ge c\},\qquad 0\le c\le1.
\]
For $d\ge2$, every level set of $|\omega_1|$ has surface measure zero.
Thus the normalized measure of $A_c$ varies continuously from one
to zero. When $\theta=0$, take $A=\varnothing$. Define
\[
F_A=\{\rho\omega:\rho\ge a,\ \omega\in A\},
\qquad \Om=\R^d\setminus F_A.
\]
On the closed set $\{|y|\ge a\}$ the map $y\mapsto y/|y|$ is
continuous, and $F_A$ is the inverse image of $A$ under this map.
Hence $F_A$ is closed and $\Om$ is open and contains $B_a(0)$.
By antipodal symmetry, the line through $0$ in direction $\omega$
meets $F_A$ for the first time at distance $a$ when $\omega\in A$.
For $\omega\notin A$ the entire line lies in $\Om$. Therefore
\[
\delta_\omega(0)=
\begin{cases}
a, & \omega\in A,\\
\infty, & \omega\notin A.
\end{cases}
\]
Within $B_r(0)$, the radial intervals in $\Om$ have lengths $a$
on $A$ and $r$ off $A$. Their radial volume fractions are $(a/r)^d$
and $1$. Averaging gives
\begin{equation}
\psi_r(0)=\theta(a/r)^d+1-\theta,
\qquad
\mathcal W_{\alpha,\ell}(0)=\theta(a+\ell)^{-\alpha}.
\label{eq:antipodal-values}
\end{equation}
Here the infinite distances off $A$ contribute zero to the weight.

For $R>r$ let $\Om_R=\Om\cap B_R(0)$. This truncation preserves the
intersection with $B_r(0)$. The directional distance remains $a$ on
$A$ and becomes $R$ off $A$. Thus the value of
$\mathcal W_{\alpha,\ell}(0)$ on $\Om_R$ is
\begin{equation}
\theta(a+\ell)^{-\alpha}+(1-\theta)(R+\ell)^{-\alpha}.
\label{eq:antipodal-truncation}
\end{equation}
For fixed $a$ and $\theta$ this converges as $R\to\infty$ to the
weight in \eqref{eq:antipodal-values} at the same volume fraction.

\begin{theorem}
\label{thm:sharpness}
Let $d\ge2$, $\alpha>0$, $r>0$, and $q\in(0,1]$. There exist an
open set $\Om\subset\R^d$ and $x\in\Om$ with $\psi_r(x)=q$ and
\begin{equation}
\frac{1}{|\Sph^{d-1}|}\int_{\Sph^{d-1}}
\delta_\omega(x)^{-\alpha}\,\dd\omega
=
\frac1{r^\alpha}\Phi_{\alpha,d}(q).
\label{eq:sharp-equality}
\end{equation}
Thus $\Phi_{\alpha,d}$ is the optimal function of $\psi_r(x)$ in
Theorem~\ref{thm:geometric}. The right side of \eqref{eq:sharp-equality}
is also the infimum over bounded open sets and points with
$\psi_r(x)=q$.
\end{theorem}

\begin{proof}
In each case we take the origin as the center and realize the
directional distances prescribed by Lemma~\ref{lem:convex-envelope}.
Recall that $q_*=\alpha/(d+\alpha)$.

If $0<q\le q_*$, a ball realizes the constant directional distance.
Take
\[
\Om=B_{r q^{1/d}}(0),
\qquad
x=0.
\]
Then
\[
\psi_r(0)=q,
\qquad
\delta_\omega(0)=r q^{1/d}
\quad\text{for every }\omega,
\]
and therefore
\[
\frac{1}{|\Sph^{d-1}|}\int_{\Sph^{d-1}}
\delta_\omega(0)^{-\alpha}\,\dd\omega
=
r^{-\alpha}q^{-\alpha/d}
=
r^{-\alpha}\Phi_{\alpha,d}(q).
\]

If $q_*<q\le1$, use the construction above with
\[
a=rq_*^{1/d},
\qquad
\theta=\frac{1-q}{1-q_*}\in[0,1).
\]
The first identity in \eqref{eq:antipodal-values} gives
\[
\psi_r(0)=\theta q_*+(1-\theta)=1-\theta(1-q_*)=q.
\]
The second identity with $\ell=0$ gives
\[
\mathcal W_{\alpha,0}(0)
=\theta a^{-\alpha}
=\frac{1-q}{1-q_*}\,q_*^{-\alpha/d}\,r^{-\alpha}
=r^{-\alpha}K_{\alpha,d}(1-q)
=r^{-\alpha}\Phi_{\alpha,d}(q).
\]
This proves equality in \eqref{eq:sharp-equality}.

The ball for $0<q\le q_*$ is bounded. For $q_*<q\le1$,
intersecting the constructed set with $B_R(0)$ preserves the volume
fraction and adds $(1-\theta)R^{-\alpha}$ to the weight by
\eqref{eq:antipodal-truncation} with $\ell=0$. Letting $R\to\infty$
gives the asserted infimum over bounded open sets, by
Theorem~\ref{thm:geometric}.
\end{proof}

\subsection{Proof of Theorem~\ref{thm:shifted-geometric}}

\begin{proof}
Fix $r>0$. The proof uses the same radial volume estimate as in
the unshifted case. The shift changes the distance at which the
ratio of the volume deficit to the weight is largest. In a direction
with boundary distance $t<r$, the missing fraction of radial volume
up to $r$ is at most $1-(t/r)^d$.
The corresponding shifted weight is $(t+\ell)^{-\alpha}$. To find the
smallest constant that bounds their ratio from above, define
\[
g_r(t)=\left(1-\left(\frac{t}{r}\right)^d\right)(t+\ell)^\alpha,
\qquad 0\le t\le r.
\]
For $0<t<r$, set $z=1-(t/r)^d$. The reciprocal $g_r(t)^{-1}$
is the slope of the line joining the origin to
$(z,(t+\ell)^{-\alpha})$. Maximizing $g_r$ thus determines the
shifted counterpart of the supporting line in the proof of
Lemma~\ref{lem:convex-envelope}. Differentiating gives
\begin{align*}
g_r'(t)
&=-\frac{d\,t^{d-1}}{r^d}(t+\ell)^\alpha
  +\alpha\left(1-\left(\frac{t}{r}\right)^d\right)(t+\ell)^{\alpha-1}\\
&=\frac{(t+\ell)^{\alpha-1}}{r^d}
  \left(\alpha r^d-(d+\alpha)t^d-d\ell t^{d-1}\right).
\end{align*}
The factor outside the parentheses is positive. The quantity in
parentheses decreases strictly and vanishes at $T_r$ by
\eqref{eq:shifted-tangency}.
Consequently, $g_r$ increases on $[0,T_r]$ and decreases on $[T_r,r]$,
so its maximum is
\[
g_r(T_r)=(1-q_r^*)(T_r+\ell)^\alpha
=C_{\alpha,d,\ell}(r)^{-1}.
\]
With $a_+:=\max\{a,0\}$, this gives
\begin{equation}
\left(1-\left(\frac{t}{r}\right)^d\right)_+
\le g_r(T_r)(t+\ell)^{-\alpha}
\qquad (0<t\le\infty).
\label{eq:shifted-scalar-linear}
\end{equation}
For $t<r$, divide $g_r(t)\le g_r(T_r)$ by $(t+\ell)^\alpha$.
Since $g_r$ has a unique maximizer, equality in this range holds
only at $t=T_r$. For finite $t\ge r$ the deficit is zero and the
right side is positive. At $t=\infty$ both sides vanish.

Subtracting the radial volume estimate \eqref{eq:polar-occupation}
from one gives
\[
1-\psi_r(x)
\le \frac{1}{\abs{\Sph^{d-1}}}\int_{\Sph^{d-1}}
\left(1-\left(\frac{\delta_\omega(x)}{r}\right)^d\right)_+
\,\dd\omega.
\]
Applying \eqref{eq:shifted-scalar-linear} with $t=\delta_\omega(x)$
bounds this average by $g_r(T_r)\mathcal W_{\alpha,\ell}(x)$.
Dividing by $g_r(T_r)$ gives
\begin{equation}
\mathcal W_{\alpha,\ell}(x)
\ge C_{\alpha,d,\ell}(r)(1-\psi_r(x)).
\label{eq:shifted-linear}
\end{equation}
This is \eqref{eq:shifted-geom} for $\psi_r(x)\ge q_r^*$.

Suppose $0<q=\psi_r(x)<q_r^*$ and put $t=rq^{1/d}$. Then $0<t<T_r$. The
first branch has value $(t+\ell)^{-\alpha}$. To obtain this value from
the linear estimate, choose a radius $s$ with $T_s=t$. Solving
\eqref{eq:shifted-tangency} at that distance gives
\begin{equation}
s=\left(\frac{(d+\alpha)t^d+d\ell t^{d-1}}{\alpha}\right)^{1/d}.
\label{eq:shifted-radius-choice}
\end{equation}
The same formula with $T_r$ in place of $t$ gives $r$. Since the
right-hand side of \eqref{eq:shifted-radius-choice} is strictly increasing
as a function of $t$, we have $0<s<r$. Its numerator is at least
$(d+\alpha)t^d>\alpha t^d$, so $t<s$. Thus $t$ lies in $(0,s)$ and
satisfies \eqref{eq:shifted-tangency} at radius $s$. Uniqueness gives
$T_s=t$ and $q_s^*=(t/s)^d$. Since $B_s(x)\subset B_r(x)$, the volume
of $\Om\cap B_s(x)$ is at most that of $\Om\cap B_r(x)$.
Dividing by the volume of the unit ball gives
\[
s^d\psi_s(x)\le r^d\psi_r(x)=r^d q=t^d,
\qquad
\psi_s(x)\le q_s^*.
\]
The deficit at $s$ is at least $1-q_s^*$. Applying
\eqref{eq:shifted-linear} at that radius, we obtain
\begin{align*}
\mathcal W_{\alpha,\ell}(x)
&\ge C_{\alpha,d,\ell}(s)(1-\psi_s(x))\\
&\ge C_{\alpha,d,\ell}(s)(1-q_s^*)\\
&=(T_s+\ell)^{-\alpha}
=(rq^{1/d}+\ell)^{-\alpha}.
\end{align*}
This is the first branch of \eqref{eq:shifted-profile}.

To attain the bound for $0<q\le q_r^*$, take $\Om=B_{r q^{1/d}}(0)$
and $x=0$, so that the distance is $rq^{1/d}$ in every direction.
The ball lies in $B_r(0)$, so
\[
\psi_r(0)=q,
\qquad
\delta_\omega(0)=r q^{1/d}\quad\text{for every }\omega\in\Sph^{d-1}.
\]
Consequently,
\[
\mathcal W_{\alpha,\ell}(0)
=(r q^{1/d}+\ell)^{-\alpha}=G_{\alpha,d,\ell,r}(q).
\]

For $q_r^*<q\le1$, we use the two equality cases of
\eqref{eq:shifted-scalar-linear}. In the construction of
\eqref{eq:antipodal-values} take $a=T_r$. Since $(T_r/r)^d=q_r^*$,
the required fraction of finite distances is the solution of
$q=\theta q_r^*+(1-\theta)$, namely
\[
\theta=\frac{1-q}{1-q_r^*}\in[0,1).
\]
This gives $\psi_r(0)=q$, and the second identity in
\eqref{eq:antipodal-values} yields
\[
\mathcal W_{\alpha,\ell}(0)
=\theta(T_r+\ell)^{-\alpha}
=C_{\alpha,d,\ell}(r)(1-q)
=G_{\alpha,d,\ell,r}(q).
\]

The ball for $0<q\le q_r^*$ is bounded. For $q_r^*<q\le1$, the
bounded sets from \eqref{eq:antipodal-truncation} preserve
$\psi_r(0)=q$ and have weight
\[
G_{\alpha,d,\ell,r}(q)+(1-\theta)(R+\ell)^{-\alpha}.
\]
The added term tends to zero as $R\to\infty$. Together with
\eqref{eq:shifted-geom}, this gives the claimed infimum over bounded
open sets. If a bounded set attained it for $q>q_r^*$, both
inequalities leading to \eqref{eq:shifted-linear} would be equalities.
In particular, \eqref{eq:shifted-scalar-linear} would be an equality
for almost every $\omega$. Since all directional distances are finite,
this would force $\delta_\omega(x)=T_r$ for almost every $\omega$.
The weight would then be $(T_r+\ell)^{-\alpha}=G_{\alpha,d,\ell,r}(q_r^*)$,
which is strictly larger than $G_{\alpha,d,\ell,r}(q)$.
\end{proof}

\section{Volume fractions at varying radii}
\label{sec:shifted}

\subsection{Proof of Theorem~\ref{thm:all-radii-sharp}}
\label{subsec:all-radii}

\begin{proof}
Let $\mu$ be normalized surface measure and put
$a_x(s)=\mu\{\omega:x+s\omega\notin\Om\}$ for $s>0$.
Polar coordinates give
\begin{equation}
D_x(r)=d\int_0^r s^{d-1}a_x(s)\,\dd s .
\label{eq:deficit-polar}
\end{equation}
If $0<s<t$ and $x+s\omega\notin\Om$, then
$\delta_\omega(x)\le s<t$. Thus $a_x(s)$ is at most the measure of
$\{\omega:\delta_\omega(x)<t\}$. By \eqref{eq:deficit-polar}, each
fraction in \eqref{eq:annular-deficit} is the weighted average
\[
\frac{\int_a^b s^{d-1}a_x(s)\,\dd s}{\int_a^b s^{d-1}\,\dd s}.
\]
Taking the supremum of these averages proves \eqref{eq:annular-tail}
for $\Om'=\Om$. The same calculation proves it for any set $\Om'$
with the same volume fractions, which determine $H_x$.
The function $H_x$ is nondecreasing by definition, lies in $[0,1]$,
and vanishes for all sufficiently small $t$ because $x$ is an interior
point. In particular the integral in \eqref{eq:all-radii-sharp} is
well defined and finite. For a finite exit distance $\delta$,
$(\delta+\ell)^{-\alpha}=\alpha\int_\delta^\infty
(t+\ell)^{-\alpha-1}\,\dd t$, with both sides zero for
$\delta=\infty$. For any open set $\Om'$ with the same profile,
Tonelli's theorem and \eqref{eq:annular-tail} give the lower bound
\[
\mathcal W_{\alpha,\ell}^{\Om'}(x)
=\alpha\int_0^\infty
\mu\{\omega:\delta_\omega^{\Om'}(x)<t\}(t+\ell)^{-\alpha-1}\,\dd t
\ge\mathcal L_{\alpha,\ell}(\Om,x).
\]

To construct an attaining set, translate $x$ to the origin. On
almost every sphere centered at the origin we rearrange the exterior
directions into a pair of antipodal caps of the same total measure.
We use a common axis at all radii. The cap pairs are then nested
by size, so the measure of their union up to radius $t$ is the
supremum of their measures over $0<s<t$.

We first express $H_0(t)$ in terms of the exterior fractions on
spheres. By \eqref{eq:deficit-polar} and the Lebesgue differentiation theorem,
\begin{equation}
H_0(t)=\operatorname*{ess\,sup}_{0<s<t} a_0(s).
\label{eq:ess-annular-deficit}
\end{equation}
Indeed, every weighted average is at most the essential supremum.
At almost every $s\in(0,t)$ the averages over intervals shrinking to
$s$ converge to $a_0(s)$, since the weight $d s^{d-1}$ is continuous
and positive there.

We need a representative of the spherical fraction that gives a
closed complement after rearrangement. The function $a_0$ is upper
semicontinuous. If $s_n\to s$, closedness of $\R^d\setminus\Om$
gives $\limsup_n\mathbf 1_{\R^d\setminus\Om}(s_n\omega)
\le \mathbf 1_{\R^d\setminus\Om}(s\omega)$.
The indicators are bounded by one, so the reverse Fatou inequality
gives $\limsup_n a_0(s_n)\le a_0(s)$. Define
\[
a^*(s)=\lim_{\varepsilon\downarrow0}
\operatorname*{ess\,sup}_{\substack{u>0\\|u-s|<\varepsilon}} a_0(u).
\]
If $a^*(s)<c$, choose a neighborhood of $s$ in which the essential
supremum of $a_0$ is below $c$. Every point in a smaller
neighborhood then has $a^*<c$, which proves upper semicontinuity.
Upper semicontinuity of $a_0$ gives $a^*(s)\le a_0(s)$ for every
$s>0$. At a Lebesgue point of $a_0$, the essential supremum over
each shrinking interval is at least its average, which converges
to $a_0(s)$. Hence $a^*=a_0$ almost everywhere. This distinction
matters for the minimization, since removing a sphere can change
$a_0$ at a single radius without changing any ball volume.

Let $F$ be the distribution function of $|\omega_1|$ under $\mu$ and
set $U(\omega)=1-F(|\omega_1|)$. For $d\ge2$, level sets of $|\omega_1|$
have zero surface measure, and every open subinterval of $(0,1)$
has positive probability. Hence $F$ is continuous and strictly
increasing on $[0,1]$. The function $U$ is continuous and even.
For $0<q<1$, the sublevel set $\{U\le q\}$ is the pair of antipodal
caps $|\omega_1|\ge F^{-1}(1-q)$ and has measure $q$ by the
definition of $F$. The same measure formula holds for $q=0,1$.
Gambicchia~\cite{Gambicchia2026} studied perimeter inequalities for
the analogous double spherical cap rearrangement of planar sets.

Choose $s_0>0$ with $B_{2s_0}(0)\subset\Om$ and put
\[
\Om_*=
\R^d\setminus\{s\omega:s\ge s_0,\ U(\omega)\le a^*(s)\}.
\]
The removed set is closed. In fact, if $s_n\omega_n$ converges with
$s_n\ge s_0$, the limit has polar coordinates $(s,\omega)$ with
$s\ge s_0$, and continuity of $U$ and upper semicontinuity of $a^*$
give $U(\omega)\le a^*(s)$ whenever $U(\omega_n)\le a^*(s_n)$.
Thus $\Om_*$ is open, contains $B_{s_0}(0)$ and is centrally
symmetric. At each radius $s\ge s_0$ the fraction removed is
$a^*(s)$. Below $s_0$, both the original and the new set contain
the entire sphere. Since $a^*=a_0$ almost everywhere, polar coordinates
give the same ball volume fractions for $\Om_*$ at all radii. If
$\Om$ is bounded, $a^*=1$ for all sufficiently large $s$, so $\Om_*$
is bounded.

It remains to compute how many lines through the origin leave $\Om_*$
before distance $t$. If $t\le s_0$ this measure is zero, and so is
$H_0(t)$. If $t>s_0$, then
\[
\sup_{s_0\le s<t}a^*(s)=H_0(t).
\]
For the upper bound, each $s\in[s_0,t)$ has a neighborhood contained
in $(0,t)$, so $a^*(s)\le H_0(t)$ by
\eqref{eq:ess-annular-deficit}. For the lower bound, use
$a^*=a_0$ almost everywhere and $a_0=0$ on $(0,s_0)$.
Because $\Om_*$ is centrally symmetric, a line in direction
$\omega$ reaches the complement before distance $t$ exactly when
$U(\omega)\le a^*(s)$ for some $s\in[s_0,t)$. The removed caps are
nested by inclusion because they are sublevel sets of the same
function $U$. Their union over $s\in[s_0,t)$ contains
$\{U<H_0(t)\}$ and is contained in $\{U\le H_0(t)\}$.
Since $\mu\{U\le q\}=q$ for $q\in[0,1]$ and $U$ has no atoms,
both sets have measure $H_0(t)$. Thus equality holds in
\eqref{eq:annular-tail} for $\Om_*$ at every $t$.
Integrating gives the claimed minimum for every $\alpha$ and $\ell$.
Translating back places the center at $x$.
\end{proof}

\subsection{Optimization over the radius}

At a small volume fraction, the linear bound at a smaller radius
can recover the first branch at the given radius. In the unshifted
case the radius change has an explicit form and applies also when
$d=1$. We first record the continuity of $r\mapsto\Psi_r$.

For $0<s\le r$, set $a=(s/r)^d$. The volume
of $\Om$ in the annulus $B_r(x)\setminus B_s(x)$ lies between zero
and the volume of the annulus. Dividing by $|B_r(x)|$ gives
\[
a\psi_s(x)\le\psi_r(x)\le a\psi_s(x)+(1-a).
\]
Since $0\le\psi_s(x)\le1$, this gives
$|\psi_r(x)-\psi_s(x)|\le1-a$ uniformly over $x\in\Om$. In particular,
\begin{equation}
|\Psi_r-\Psi_s|\le\sup_{x\in\Om}|\psi_r(x)-\psi_s(x)|\le1-a.
\label{eq:volume-continuity}
\end{equation}
It follows that $r\mapsto\Psi_r$ is continuous on $(0,\infty)$.

\begin{proposition}
\label{prop:radius-optimization}
Let $d\in\mathbb N$, let $\Om\subset\R^d$ be nonempty and open, and let
$\alpha>0$. With $K_{\alpha,d}$ as in
\eqref{eq:K-alpha-d}, one has
\begin{equation}
\sup_{r>0}\frac{\Phi_{\alpha,d}(\Psi_r)}{r^\alpha}
=
K_{\alpha,d}\sup_{r>0}\frac{1-\Psi_r}{r^\alpha}.
\label{eq:radius-optimization}
\end{equation}
\end{proposition}

\begin{proof}
For $0<s\le r$, the inclusion $B_s(x)\subset B_r(x)$ gives
$s^d\psi_s(x)\le r^d\psi_r(x)$, since both sides are intersection
volumes divided by the volume of the unit ball. Taking the supremum
over $x\in\Om$ yields
\begin{equation}
\Psi_s\le\left(\frac rs\right)^d\Psi_r.
\label{eq:density-scale}
\end{equation}
For $q<q_*$, the line $K_{\alpha,d}(1-q)$ lies below the convex function
$q^{-\alpha/d}$ by tangency at $q_*$. For $q\ge q_*$, the line equals
$\Phi_{\alpha,d}(q)$. Thus $\Phi_{\alpha,d}(q)\ge K_{\alpha,d}(1-q)$
for all $q\in(0,1]$. Applying this at each radius and taking suprema
gives one direction of \eqref{eq:radius-optimization}.

For the reverse inequality, fix $r>0$ and write $q=\Psi_r>0$.

If $q\ge q_*$, the two expressions already agree at radius $r$,
\[
\frac{\Phi_{\alpha,d}(\Psi_r)}{r^\alpha}
=K_{\alpha,d}\frac{1-\Psi_r}{r^\alpha}.
\]

If $q<q_*$, choose $s<r$ so that the upper bound on $\Psi_s$ in
\eqref{eq:density-scale} equals the transition value $q_*$. The
equation $(r/s)^d q=q_*$ gives $s=r(q/q_*)^{1/d}$. Since
$\Psi_s\le q_*$, the linear expression at radius $s$ satisfies
\begin{align*}
K_{\alpha,d}\frac{1-\Psi_s}{s^\alpha}
&\ge K_{\alpha,d}\frac{1-q_*}{s^\alpha}\\
&=\frac{q^{-\alpha/d}}{r^\alpha}
=\frac{\Phi_{\alpha,d}(\Psi_r)}{r^\alpha}.
\end{align*}
The first equality uses the formula for $s$ and the identity
$K_{\alpha,d}(1-q_*)=q_*^{-\alpha/d}$, which follows from agreement
of the two branches at $q_*$. Thus for each $r$, the sharp bound at $r$
is bounded above by a value of the linear expression. Taking the
supremum over $r>0$ proves the reverse inequality.
\end{proof}

For the shifted bound, the choice of the smaller radius is determined
by the distance $T_r$. The same argument applies to volume fractions
at a fixed center as well as to their suprema over centers.

\begin{proposition}
\label{prop:shifted-radius-optimization}
Let $d\ge2$, $\alpha>0$, $\ell\ge0$, and let $\Om\subset\R^d$ be
nonempty and open. Then
\begin{equation}
\sup_{r>0}G_{\alpha,d,\ell,r}(\Psi_r)
=\sup_{r>0}C_{\alpha,d,\ell}(r)(1-\Psi_r).
\label{eq:shifted-radius-optimization}
\end{equation}
The same identity holds with $\Psi_r$ replaced by $\psi_r(x)$ at any
fixed $x\in\Om$.
\end{proposition}

\begin{proof}
Let $M_r$ denote either $\Psi_r$ or $\psi_r(x)$ for a fixed $x\in\Om$.
In both cases, inclusion of the balls gives $s^dM_s\le r^dM_r$ for
$0<s\le r$, as in the proof of \eqref{eq:density-scale}.

For $q\ge q_r^*$, $G_{\alpha,d,\ell,r}(q)$ equals the linear expression.
For $q<q_r^*$, put $t=rq^{1/d}$. The maximum of $g_r$ computed in the
proof of Theorem~\ref{thm:shifted-geometric} gives
\[
(1-q)(t+\ell)^\alpha=g_r(t)\le g_r(T_r)
=C_{\alpha,d,\ell}(r)^{-1}.
\]
Rearranging yields
\[
C_{\alpha,d,\ell}(r)(1-q)\le(t+\ell)^{-\alpha}
=G_{\alpha,d,\ell,r}(q).
\]
Thus $G_{\alpha,d,\ell,r}(q)$ is at least the linear expression for every $q\in(0,1]$.
Applying this at $q=M_r$ and taking suprema gives
\[
\sup_{r>0}G_{\alpha,d,\ell,r}(M_r)
\ge \sup_{r>0}C_{\alpha,d,\ell}(r)(1-M_r).
\]

For the reverse inequality, fix $r>0$ and write $q=M_r$. If $q\ge q_r^*$,
$G_{\alpha,d,\ell,r}(q)$ equals the linear expression at $r$. If $q<q_r^*$, set
$t=rq^{1/d}$ and choose $s$ by \eqref{eq:shifted-radius-choice}, so that
$s<r$ and $T_s=t$. The scaling bound then gives
\[
M_s\le (r/s)^dM_r=t^d/s^d=q_s^*.
\]
Consequently, the linear expression at $s$ satisfies
\begin{align*}
C_{\alpha,d,\ell}(s)(1-M_s)
&\ge C_{\alpha,d,\ell}(s)(1-q_s^*)\\
&=(T_s+\ell)^{-\alpha}\\
&=G_{\alpha,d,\ell,r}(M_r).
\end{align*}
In either case, $G_{\alpha,d,\ell,r}(M_r)$ is bounded by the supremum of
the linear expression. Taking the supremum over $r$ proves the
reverse inequality for both choices of $M_r$.
\end{proof}

\subsection{An integral bound in terms of ball volumes}

The integral bound below uses the deficit at each radius separately.
Its positive kernel allows the use of $\Psi_r$, the supremum of
the ball volume fractions over centers at each radius, in the
spectral estimates.

\begin{proposition}
\label{prop:all-radii}
Let $d\ge2$, $\alpha>0$ and $\ell\ge0$, and let $\Om\subset\R^d$ be
nonempty and open. Put
\begin{equation}
k_{\alpha,d,\ell}(r)
=\frac{\alpha\big((d+\alpha)r+(d-1)\ell\big)}
       {d(r+\ell)^{\alpha+2}},\qquad r>0.
\label{eq:radial-kernel}
\end{equation}
Then, for every $x\in\Om$,
\begin{equation}
\mathcal W_{\alpha,\ell}(x)\ge
\int_0^\infty k_{\alpha,d,\ell}(r)(1-\psi_r(x))\,\dd r .
\label{eq:all-radii-pointwise}
\end{equation}
In particular,
\begin{equation}
\begin{aligned}
\inf_{x\in\Om}\mathcal W_{\alpha,\ell}(x)
&\ge \inf_{x\in\Om}\mathcal L_{\alpha,\ell}(\Om,x)\\
&\ge \max\left\{
\sup_{r>0}G_{\alpha,d,\ell,r}(\Psi_r),\,
\int_0^\infty k_{\alpha,d,\ell}(r)(1-\Psi_r)\,\dd r
\right\}.
\end{aligned}
\label{eq:all-radii-uniform}
\end{equation}
At zero shift the kernel is
$k_{\alpha,d,0}(r)=\alpha(d+\alpha)/(d\,r^{\alpha+1})$.
\end{proposition}

\begin{proof}
Fix $x\in\Om$ and let $a_x(r)=\mu\{\omega:x+r\omega\notin\Om\}$,
where $\mu$ is normalized surface measure. If $x+r\omega\notin\Om$,
then $\delta_\omega(x)\le r$. Writing each inverse power as an
integral from $\delta_\omega(x)$ to infinity and applying Tonelli's
theorem gives
\begin{equation}
\mathcal W_{\alpha,\ell}(x)
=\int_0^\infty \frac{\alpha\,\mu\{\omega:\delta_\omega(x)\le r\}}
                        {(r+\ell)^{\alpha+1}}\,\dd r
\ge\int_0^\infty \frac{\alpha\,a_x(r)}{(r+\ell)^{\alpha+1}}\,\dd r .
\label{eq:spherical-deficit}
\end{equation}
By \eqref{eq:deficit-polar}, $D_x$ is locally absolutely continuous
and $D_x'(r)=d r^{d-1}a_x(r)$ almost everywhere. Put
$f(r)=r^{1-d}(r+\ell)^{-\alpha-1}$. For $0<\varepsilon<R$,
integration by parts gives
\[
\int_\varepsilon^R \frac{\alpha\,a_x(r)}{(r+\ell)^{\alpha+1}}\,\dd r
=\frac{\alpha}{d}[f(r)D_x(r)]_\varepsilon^R
-\frac{\alpha}{d}\int_\varepsilon^R D_x(r)f'(r)\,\dd r.
\]
The lower boundary term is zero for small enough $\varepsilon$,
since $B_\varepsilon(x)\subset\Om$. At the upper endpoint,
$0\le f(R)D_x(R)\le R(R+\ell)^{-\alpha-1}\to0$.
The product rule gives
\[
-\frac{\alpha}{d}r^d f'(r)
=\frac{\alpha}{d}\left(
\frac{d-1}{(r+\ell)^{\alpha+1}}
+\frac{(\alpha+1)r}{(r+\ell)^{\alpha+2}}\right)
=k_{\alpha,d,\ell}(r).
\]
Since $D_x(r)=r^d(1-\psi_r(x))$ and the kernel is positive,
monotone convergence lets $\varepsilon\downarrow0$ and $R\uparrow\infty$.
Using the resulting identity in \eqref{eq:spherical-deficit} yields
\eqref{eq:all-radii-pointwise}.

The first inequality in \eqref{eq:all-radii-uniform} follows from
Theorem~\ref{thm:all-radii-sharp}. To prove the second, fix $x\in\Om$
and apply \eqref{eq:shifted-geom} and \eqref{eq:all-radii-pointwise} to
an attaining set from that theorem. Its weight at $x$ is
$\mathcal L_{\alpha,\ell}(\Om,x)$ and its volume fractions are
$\psi_r(x)$ at every radius. Thus $\mathcal L_{\alpha,\ell}(\Om,x)$ is
at least $\sup_{r>0}G_{\alpha,d,\ell,r}(\psi_r(x))$ and at least the
integral in \eqref{eq:all-radii-pointwise}. Since $\psi_r(x)\le\Psi_r$,
$G_{\alpha,d,\ell,r}$ is decreasing, and the kernel is positive,
replacing $\psi_r(x)$ by $\Psi_r$ gives the two lower bounds in
\eqref{eq:all-radii-uniform} independently of $x$. Taking the infimum
over $x$ completes the proof.
\end{proof}

At a fixed center, $\mathcal L_{\alpha,\ell}(\Om,x)$ can strictly
exceed both the integral and the optimized fixed-radius bound.

\begin{example}
\label{ex:annulus-all-radii}
Let $d=\alpha=2$, $\ell=0$ and
$\Om=B_1(0)\cup\{y:\sqrt{3/2}<|y|<2\}$. At the origin, $H_0(t)=0$
for $t\le1$. For $t>1$, choose $1<b<\min\{t,\sqrt{3/2}\}$.
The annulus with radii $1$ and $b$ lies outside $\Om$, so $H_0(t)=1$.
Thus $\mathcal L_{2,0}(\Om,0)=1$.

The spherical deficit $a_0(r)$ is one on $(1,\sqrt{3/2})$ and
$(2,\infty)$ and zero elsewhere, up to the endpoints. The
integration-by-parts identity in the proof of
Proposition~\ref{prop:all-radii} gives
\begin{align*}
4\int_0^\infty (1-\psi_r(0))\,\frac{\dd r}{r^3}
&=2\int_0^\infty a_0(r)\,\frac{\dd r}{r^3}\\
&=2\int_1^{\sqrt{3/2}}\frac{\dd r}{r^3}
  +2\int_2^\infty\frac{\dd r}{r^3}=\frac7{12}.
\end{align*}
By Proposition~\ref{prop:shifted-radius-optimization}, the optimized
single-radius bound at the origin is the supremum of
$4(1-\psi_r(0))/r^2$. With $s=r^2$, the areas missing from $B_r(0)$,
divided by $\pi$, are respectively $0$, $s-1$, $1/2$ and $s-7/2$.
Thus the expression to maximize is
\[
\begin{cases}
0, & 0<s\le1,\\
4(s-1)/s^2, & 1\le s\le3/2,\\
2/s^2, & 3/2\le s\le4,\\
4(s-7/2)/s^2, & s\ge4.
\end{cases}
\]
The second branch has derivative $4(2-s)/s^3$, so its maximum is
$8/9$ at $s=3/2$. The third decreases from that value. The last
has derivative $4(7-s)/s^3$ and reaches $2/7$ at $s=7$. The
supremum is therefore $8/9$. The rays leave the inner ball at
distance one and later enter the outer annulus. Thus the exact
minimum for the full volume profile satisfies
$\mathcal L_{2,0}(\Om,0)=1>\max\{8/9,7/12\}$.
\end{example}

In a convex domain a ray stays outside once it has left the domain.
If the domain is also centrally symmetric, its center of symmetry
maximizes the ball fraction at every radius. Together these facts
give equality in the integral bound.

\begin{proposition}
\label{prop:symmetric-convex}
Let $d\ge2$, $\alpha>0$ and $\ell\ge0$. If $\Om\subset\R^d$ is a
nonempty bounded centrally symmetric open convex set, then
\begin{equation}
\inf_{x\in\Om}\mathcal W_{\alpha,\ell}(x)
=\int_0^\infty k_{\alpha,d,\ell}(r)(1-\Psi_r)\,\dd r .
\label{eq:symmetric-convex}
\end{equation}
If $\Om$ is not a ball, the integral is strictly larger than
$\sup_{r>0}G_{\alpha,d,\ell,r}(\Psi_r)$.
\end{proposition}

\begin{proof}
Translate the center of symmetry to the origin. For $r>0$ and
$x\in\R^d$, put $A_x=\Om\cap B_r(x)$. If $A_x$ is nonempty,
central symmetry gives $A_{-x}=-A_x$. Given $y\in A_x$ and
$z\in A_{-x}$, their midpoint lies in $\Om$ by convexity and
in $B_r(0)$ because $|y+z|\le|y-x|+|z+x|<2r$.
The Brunn--Minkowski inequality then gives
\[
|A_0|^{1/d}\ge
\left|\frac{A_x+A_{-x}}2\right|^{1/d}
\ge\frac{|A_x|^{1/d}+|A_{-x}|^{1/d}}2=|A_x|^{1/d}.
\]
The same conclusion holds when $A_x$ is empty. In particular,
$\Psi_r=\psi_r(0)$ for every $r$.

Let $\rho(\omega)$ be the radial distance from the origin to
$\partial\Om$ in direction $\omega$. Convexity implies that
$r\omega\notin\Om$ exactly when $r\ge\rho(\omega)$. By symmetry both
exit distances equal $\rho(\omega)$. Therefore
\eqref{eq:spherical-deficit} is an equality at the origin. Since
$\Psi_r=\psi_r(0)$, the integral in \eqref{eq:symmetric-convex} equals
$\mathcal W_{\alpha,\ell}(0)$ and bounds the weight from below at
every center. This proves \eqref{eq:symmetric-convex}. At zero shift
the spherical average of $\rho^{-\alpha}$ is the radial moment in
O'Hara's formula for the renormalized Riesz potential~\cite{OHara2012}.

For the strict comparison, the radial distance is continuous.
Indeed, to see subadditivity of the Minkowski functional of $\Om$,
take $y\in s\Om$ and $z\in t\Om$. Convexity gives
$y+z\in(s+t)\Om$. It is bounded above by $|y|/a$ when
$B_a(0)\subset\Om$. Applying subadditivity to $y=z+(y-z)$ and
$z=y+(z-y)$ gives Lipschitz continuity. On the unit sphere the
functional is positive by boundedness and equals $1/\rho$. For $r>0$, polar
coordinates give
\[
1-\Psi_r=\int_{\Sph^{d-1}}\big(1-(\rho(\omega)/r)^d\big)_+\,\dd\mu(\omega).
\]
Equality in \eqref{eq:shifted-scalar-linear} at a finite positive
distance occurs only at $T_r$. Since $\Om$ is not a ball, $\rho$ is
nonconstant. By continuity it differs from $T_r$ on a set of positive
surface measure. Integration gives
$C_{\alpha,d,\ell}(r)(1-\Psi_r)<\mathcal W_{\alpha,\ell}(0)$.

The left-hand side is continuous in $r$ by
\eqref{eq:shifted-tangency}, \eqref{eq:shifted-linear-coef} and the
continuity of $\Psi_r$. It vanishes for $r\le a$ and is positive for
large $r$ because $\Om$ is bounded. Evaluating $g_r$ at $r/2$ gives
$C_{\alpha,d,\ell}(r)\le(1-2^{-d})^{-1}(r/2+\ell)^{-\alpha}$,
so the product tends to zero at infinity. Choose a radius where
the product is positive. Outside a compact interval containing
that radius it is smaller than half this positive value. It
therefore attains its supremum at a finite positive radius,
where the comparison is strict.
Proposition~\ref{prop:shifted-radius-optimization} gives the asserted
comparison with the optimized bound.
\end{proof}

\section{Spectral applications}
\label{sec:applications}

Hardy and Hardy--Rellich inequalities bound the corresponding
quadratic forms below by weighted integrals. A uniform lower bound for
the weight gives a form bound $h[u]\ge c\|u\|_{L^2(\Om)}^2$. Dividing
by the squared norm for nonzero $u$ and taking the infimum over
the form domain gives the spectral bound $\lambda\ge c$.
When only the suprema $\Psi_r$ are available, we use
\eqref{eq:all-radii-uniform}. When the fractions are known at each
center, we use Theorem~\ref{thm:all-radii-sharp} at that center before
taking the infimum over centers.

\subsection{The Robin Laplacian}

Let $d\ge2$, let $\Om\subset\R^d$ be nonempty and open with
uniformly Lipschitz boundary, and let $\sigma>0$ be constant.
We use the Robin form \eqref{eq:robin-form}.
Let $\mu$ be normalized surface measure on $\Sph^{d-1}$ and define
\begin{equation}
\beta_d:=d\int_{\Sph^{d-1}}|\omega_1|\,\dd\mu(\omega),
\qquad
\ell_\sigma:=\frac{\beta_d}{2\sigma}.
\label{eq:robin-shift}
\end{equation}
The factor $\beta_d$ arises from averaging the boundary projection
factor. It determines the shift for the Robin form with coefficient
$\sigma$.

\begin{lemma}
\label{lem:robin-hardy}
For every $u\in H^1(\Om)$,
\begin{equation}
h_{\Om}^{\mathrm{R}\sigma}[u]
\ge
\frac d4\int_\Om |u(x)|^2
\left(\int_{\Sph^{d-1}}
\frac{\dd\mu(\omega)}{(\delta_\omega(x)+\ell_\sigma)^2}\right)\dd x.
\label{eq:robin-hardy}
\end{equation}
If $\Om$ is bounded and $D=\operatorname{diam}\Om$, the right-hand
side can be increased by
$d\|u\|_{L^2(\Om)}^2/[4(D/2+\ell_\sigma)^2]$.
\end{lemma}

\begin{proof}
We use the interval inequality of
Kova\v{r}\'{i}k and Laptev~\cite[Lemma~4.1]{KovarikLaptev2012}.
Spherical averaging gives a factor $1/d$ for the gradient term
and $\beta_d/d$ for the projected boundary term. To recover the
Robin form with boundary coefficient $\sigma$, we keep the interval
parameter $\gamma>0$ free until we average, then take
$\gamma=\sigma/\beta_d$. Put $L=(2\gamma)^{-1}$. On $[0,h]$, with $h>0$,
the function $f(t)=(t+L)^{1/2}$ satisfies $f'(0)/f(0)=\gamma$ and
$-f''(t)/f(t)=1/[4(t+L)^2]$. Given $v\in H^1(0,h)$, set
$w=v/f$. Expanding $v'=fw'+f'w$ gives
\[
|v'|^2=f^2|w'|^2+f'^2|w|^2+ff'(|w|^2)'.
\]
Since $(ff')'=f'^2+ff''$, integration by parts yields
\[
\int_0^h |v'|^2\,\dd t
=\int_0^h f^2|w'|^2\,\dd t
-\int_0^h\frac{f''}{f}|v|^2\,\dd t
+\left[\frac{f'}{f}|v|^2\right]_0^h.
\]
Substituting $f''/f$ and the endpoint values of $f'/f$ gives
\begin{align*}
&\int_0^h |v'(t)|^2\,\dd t+\gamma |v(0)|^2
-\frac14\int_0^h\frac{|v(t)|^2}{(t+L)^2}\,\dd t\\
&\qquad=
\int_0^h f(t)^2\left|\left(\frac{v}{f}\right)'(t)\right|^2\,\dd t
+\frac{|v(h)|^2}{2(h+L)}\ge0.
\end{align*}
The identity uses $f'(0)/f(0)=\gamma$ to cancel the boundary
contribution at $0$. The contribution at $h$ is nonnegative for
every trace value $v(h)$.

For a bounded interval $I=(a,b)$, apply the estimate on each half
with $h=(b-a)/2$. Measure distance from $a$ on the left half and from
$b$ on the right half, so that the two weights use $t-a$ and $b-t$,
respectively. Discard the two nonnegative midpoint contributions.
For a half-line, measure distance from its finite endpoint and take
$h$ beyond the support of $v$. Let $D_I(t)$ be the distance to the
nearest finite endpoint, with
$D_I(t)=\infty$ if $I=\R$. The interval estimate is then
\begin{equation}
\int_I |v'(t)|^2\,\dd t+\gamma\sum_{c\in\partial I}|v(c)|^2
\ge
\frac14\int_I\frac{|v(t)|^2}{(D_I(t)+L)^2}\,\dd t.
\label{eq:robin-line-hardy}
\end{equation}
For unbounded $I$, the applications below use restrictions of smooth
compactly supported functions. If $I=\R$, the right-hand side is zero.

For a finite interval, a stronger estimate follows by using
\cite[Lemma~3.8]{KovarikLaptev2012} on the two halves. On a half
of length $h$, with $t$ measured from the finite endpoint, that
lemma gives
\[
\int_0^h |v'(t)|^2\,\dd t+\gamma |v(0)|^2
\ge \frac14\int_0^h
\left(\frac1{(t+L)^2}+\frac1{(h+L)^2}\right)|v(t)|^2\,\dd t
+\frac{|v(0)|^2}{2(h+L)}.
\]
To obtain the two additional terms, put
$A(t)=(t+L)^{-1}$, $B=(h+L)^{-1}$ and $g=(A-B)/2$.
Expanding the square and integrating the term $-g(|v|^2)'$
gives
\[
0\le\int_0^h |v'-gv|^2\,\dd t
=\int_0^h |v'|^2\,\dd t+g(0)|v(0)|^2
+\int_0^h (g'+g^2)|v|^2\,\dd t ,
\]
since $g(h)=0$. Because $A\ge B$,
$-g'-g^2=(A^2+2AB-B^2)/4\ge(A^2+B^2)/4$.
Combining this with $\gamma-g(0)=B/2$ proves the displayed
bound.
Dropping the last term and adding the two half-interval inequalities
increases the right-hand side of \eqref{eq:robin-line-hardy} by
$\|v\|_{L^2(I)}^2/[4(h+L)^2]$, where $h=|I|/2$.

Take $u$ to be the restriction to $\Om$ of a function in
$C_c^\infty(\R^d)$. Fix $\omega\in\Sph^{d-1}$ and write
$x=y+t\omega$ with $y\in\omega^\perp$. The slice
\[
\Om_{\omega,y}:=\{t\in\R:y+t\omega\in\Om\}
\]
is a disjoint union of open intervals. For $t$ in a component $I$,
the first exit toward an endpoint is at that endpoint. An unbounded
end has infinite exit distance and contributes no boundary value.
Taking the nearer endpoint gives $D_I(t)=\delta_\omega(y+t\omega)$.
Apply \eqref{eq:robin-line-hardy} to $v(t)=u(y+t\omega)$ on each
component, sum over the components, and integrate in $y$.

To evaluate the endpoint terms, consider orthogonal projection
$\pi_\omega:\partial\Om\to\omega^\perp$. Its tangential Jacobian
at a differentiability point $z$ of the boundary is
$|\nu(z)\cdot\omega|$, where $\nu(z)$ is the outer unit normal.
Indeed, the projection from $\nu(z)^\perp$ onto $\omega^\perp$
fixes their common directions, and in the plane of the two
normals its length factor is $|\nu(z)\cdot\omega|$. For parallel
normals the projection is the identity on the tangent plane.
The projections of the exceptional boundary points have zero
$(d-1)$-dimensional measure. By the area formula, the same holds
for the projection of the set where $|\nu(z)\cdot\omega|=0$.
Consequently, for almost every $y$, each intersection of
$y+\R\omega$ with $\partial\Om$ is transverse. Near such an
intersection $z$, write $F(\xi)=\xi_d-\varphi(\xi')$ in local graph
coordinates, with $\Om$ on one side of $\{F=0\}$. The expansion
$F(z+s\omega)=s\nabla F(z)\cdot\omega+o(s)$ and the nonzero normal
component $\nabla F(z)\cdot\omega$ show that $F$ changes sign across
$s=0$. Each intersection therefore contributes one slice endpoint,
and the area formula gives
\[
\int_{\omega^\perp}\sum_{t\in\partial\Om_{\omega,y}}
    |u(y+t\omega)|^2\,\dd y
=\int_{\partial\Om}|\nu(z)\cdot\omega|\,|u(z)|^2\,\dd S(z).
\]
For the volume integrals, Fubini's theorem identifies integration in
$(y,t)$ with integration in $\Om$. Since
$\partial_t u(y+t\omega)=\omega\cdot\nabla u(y+t\omega)$,
we obtain
\begin{align}
&\int_\Om |\omega\cdot\nabla u(x)|^2\,\dd x
+\gamma\int_{\partial\Om}|\nu(z)\cdot\omega|\,|u(z)|^2\,\dd S(z)
\notag\\
&\qquad\ge
\frac14\int_\Om
\frac{|u(x)|^2}{(\delta_\omega(x)+L)^2}\,\dd x.
\label{eq:robin-one-direction}
\end{align}
If $\Om$ is bounded, then $|I|\le D$ on every slice component.
The additional term on each component is therefore at least
$\|v\|_{L^2(I)}^2/[4(D/2+L)^2]$. By Fubini's theorem,
\[
\int_{\omega^\perp}\sum_I\int_I |u(y+t\omega)|^2\,\dd t\,\dd y
=\|u\|_{L^2(\Om)}^2,
\]
where the sum is over the components of the slice at $y$. Thus the
extra term in \eqref{eq:robin-one-direction} is
$\|u\|_{L^2(\Om)}^2/[4(D/2+L)^2]$.

To average the directional estimate over $\omega$, observe that
the spherical averages of $\omega_i\omega_j$ vanish for $i\ne j$ by
reflection symmetry. The diagonal averages are equal and sum to one,
so each is $1/d$. Expand $|\omega\cdot\nabla u|^2$ using these moments.
For the boundary term, rotate the unit normal to the first coordinate
vector. These calculations give
\[
\int_{\Sph^{d-1}}|\omega\cdot\nabla u(x)|^2\,\dd\mu(\omega)
=\frac1d|\nabla u(x)|^2,
\qquad
\int_{\Sph^{d-1}}|\nu(z)\cdot\omega|\,\dd\mu(\omega)
=\frac{\beta_d}{d}.
\]
For fixed $t\ge0$, the condition $\delta_\omega(x)>t$ means that the
compact segment $\{x+s\omega:|s|\le t\}$ lies in $\Om$. This condition
is preserved under sufficiently small changes in $x$ and $\omega$,
so the distance is jointly lower semicontinuous and the weight is
measurable. Average \eqref{eq:robin-one-direction} with respect to
$\mu$, using Tonelli's theorem. The coefficients of the gradient
and boundary integrals on the left are $1/d$ and
$\gamma\beta_d/d$, respectively. Multiplication by $d$ gives
gradient coefficient one, boundary coefficient $\gamma\beta_d$,
and coefficient $d/4$ in front of the averaged weight. Choosing
$\gamma=\sigma/\beta_d$ gives the boundary coefficient $\sigma$
and shift $L=(2\gamma)^{-1}=\ell_\sigma$. This proves
\eqref{eq:robin-hardy} for the smooth restrictions.
In the bounded case the additional term is independent of the direction.
It is unchanged by the normalized average, and multiplication by $d$
gives the claimed term with denominator $4(D/2+\ell_\sigma)^2$.

To pass to $H^1(\Om)$, extend $u$ to $H^1(\R^d)$ by the Sobolev
extension theorem for uniformly Lipschitz domains and approximate
the extension by $C_c^\infty(\R^d)$ functions. Their restrictions
converge in $H^1(\Om)$, and the bounded trace map gives convergence
in $L^2(\partial\Om)$ of the traces. If $u_j$ denotes these
restrictions, then Cauchy--Schwarz gives
\[
\||u_j|^2-|u|^2\|_{L^1(\Om)}
\le \|u_j-u\|_{L^2(\Om)}(\|u_j\|_{L^2(\Om)}+\|u\|_{L^2(\Om)}).
\]
The Hardy weight is bounded by $\ell_\sigma^{-2}$, so this
also gives convergence of the weighted integral. The squared
$L^2$ norms converge as well. Both form estimates therefore hold
for all $u\in H^1(\Om)$.
\end{proof}

\begin{theorem}
\label{thm:robin}
Let $d\ge2$ and let $\Om\subset\R^d$ be nonempty and open with
uniformly Lipschitz boundary. If $\sigma>0$, then for every $r>0$,
\begin{equation}
\lambda_\Om^{\mathrm{R}\sigma}
\ge\frac d4 G_{2,d,\ell_\sigma,r}(\Psi_r).
\label{eq:robin-main}
\end{equation}
Optimizing over $r$ gives
\begin{equation}
\lambda_\Om^{\mathrm{R}\sigma}
\ge\frac d4\sup_{r>0}G_{2,d,\ell_\sigma,r}(\Psi_r)
=\frac d4\sup_{r>0} C_{2,d,\ell_\sigma}(r)(1-\Psi_r).
\label{eq:robin-optimized}
\end{equation}
There is also the bound
\begin{equation}
\lambda_\Om^{\mathrm{R}\sigma}
\ge \frac d4\int_0^\infty k_{2,d,\ell_\sigma}(r)(1-\Psi_r)\,\dd r .
\label{eq:robin-integral}
\end{equation}
If $\Om$ is bounded with $D=\operatorname{diam}\Om$, each of the
three lower bounds can be increased by
\begin{equation}
\frac{d}{4(D/2+\ell_\sigma)^2}.
\label{eq:robin-diameter-remainder}
\end{equation}
Using the volume fractions at each center also gives
\[
\lambda_\Om^{\mathrm{R}\sigma}
\ge \frac d4\inf_{x\in\Om}\mathcal L_{2,\ell_\sigma}(\Om,x),
\]
with the same additional term \eqref{eq:robin-diameter-remainder}
when $\Om$ is bounded.
\end{theorem}

\begin{proof}
Fix $r>0$. For each $x\in\Om$, Theorem~\ref{thm:shifted-geometric}
with $\alpha=2$ and $\ell=\ell_\sigma$, together with
$\psi_r(x)\le\Psi_r$ and the monotonicity of $G_{2,d,\ell_\sigma,r}$,
gives
\[
\int_{\Sph^{d-1}}\frac{\dd\mu(\omega)}{(\delta_\omega(x)+\ell_\sigma)^2}
\ge G_{2,d,\ell_\sigma,r}(\psi_r(x))
\ge G_{2,d,\ell_\sigma,r}(\Psi_r).
\]
Substituting the last bound into \eqref{eq:robin-hardy} gives
\[
h_\Om^{\mathrm{R}\sigma}[u]
\ge \frac d4 G_{2,d,\ell_\sigma,r}(\Psi_r)\|u\|_{L^2(\Om)}^2
\qquad (u\in H^1(\Om)).
\]
The variational principle gives \eqref{eq:robin-main}. Taking the
supremum over $r>0$ and applying
Proposition~\ref{prop:shifted-radius-optimization} gives
\eqref{eq:robin-optimized}. Using the integral in
\eqref{eq:all-radii-uniform} as the uniform estimate for the same
Hardy weight gives \eqref{eq:robin-integral}. On a bounded centrally
symmetric convex set other than a ball, this bound strictly exceeds
\eqref{eq:robin-optimized} by Proposition~\ref{prop:symmetric-convex}.

Theorem~\ref{thm:all-radii-sharp} also gives
$\mathcal W_{2,\ell_\sigma}(x)\ge\mathcal L_{2,\ell_\sigma}(\Om,x)$
at each $x$. Taking the infimum over centers yields a constant
lower bound for the Hardy weight in \eqref{eq:robin-hardy}.
Thus the Robin form is bounded below by $d/4$ times this infimum
times $\|u\|_{L^2(\Om)}^2$, and the variational principle gives the
claimed bound using the volume fractions at each center.

For bounded $\Om$, the extra term in Lemma~\ref{lem:robin-hardy}
is a constant times $\|u\|_{L^2(\Om)}^2$. Dividing by this squared
norm carries the constant into each Rayleigh quotient, so the same
term \eqref{eq:robin-diameter-remainder} can be added to the
fixed-radius and integral bounds and to the bound using the
volume fractions at each center. As the constant is independent of $r$,
it also survives taking the supremum over $r$.
\end{proof}

\subsection{Dirichlet polyharmonic operators}

Let $m\in\mathbb N$. We consider the Dirichlet realization of
$(-\Delta)^m$ in $L^2(\Om)$, with quadratic form
\begin{equation}
h_{\Om}^{(m)}[u]
=
\begin{cases}
\displaystyle\int_{\Om}\abs{(-\Delta)^{m/2}u}^2\,\dd x,
& m\ \text{even},\\[6pt]
\displaystyle\int_{\Om}\abs{\nabla(-\Delta)^{(m-1)/2}u}^2\,\dd x,
& m\ \text{odd},
\end{cases}
\end{equation}
on $H_0^m(\Om)$, the closure of $C_c^\infty(\Om)$ in $H^m(\Om)$.
Write $\lambda_{\Om}^{(m)}$ for the bottom of the spectrum. For $m=1$,
we also write $\lambda_{\Om}^{\mathrm D}=\lambda_{\Om}^{(1)}$ for the
bottom of the spectrum of the Dirichlet Laplacian.

The local-volume estimate follows from Owen's Hardy--Rellich
inequality. For $x\in\Om$, set
\begin{equation}
M_{\Om}^{(m)}(x)
:=
\frac{1}{|\Sph^{d-1}|}\int_{\Sph^{d-1}}
\delta_\omega(x)^{-2m}\,\dd\omega.
\end{equation}
Owen's Hardy--Rellich inequality~\cite{Owen1999} states that
\begin{equation}
h_{\Om}^{(m)}[u]
\ge C_{m,d}\int_{\Om}\abs{u(x)}^2M_{\Om}^{(m)}(x)\,\dd x
\qquad\text{for all }u\in H_0^m(\Om),
\label{eq:Owen}
\end{equation}
where
\begin{equation}
C_{m,d}
=
\frac{(d+2m-2)(d+2m-4)\cdots d\,(2m-1)(2m-3)\cdots1}{4^m}.
\label{eq:Cmd}
\end{equation}

\begin{theorem}
\label{thm:polyharmonic}
Let $d\ge2$, let $\Om\subset\R^d$ be nonempty and open, and let
$m\in\mathbb N$. For every $r>0$,
\begin{equation}
\lambda_{\Om}^{(m)}
\ge \frac{C_{m,d}}{r^{2m}}\Phi_{2m,d}(\Psi_r).
\label{eq:poly-main}
\end{equation}
In explicit form,
\begin{equation}
\lambda_{\Om}^{(m)}
\ge \frac{C_{m,d}}{r^{2m}}
\begin{cases}
\Psi_r^{-2m/d},
& 0<\Psi_r\le \dfrac{2m}{d+2m},\\[8pt]
c_{m,d}(1-\Psi_r),
& \dfrac{2m}{d+2m}\le\Psi_r\le1,
\end{cases}
\label{eq:poly-piecewise}
\end{equation}
where
\begin{equation}
c_{m,d}
=\frac{d+2m}{d}\left(\frac{d+2m}{2m}\right)^{2m/d}.
\label{eq:cmd}
\end{equation}
Taking the supremum over $r>0$ gives
\begin{equation}
\lambda_{\Om}^{(m)}
\ge C_{m,d}\sup_{r>0}\frac{\Phi_{2m,d}(\Psi_r)}{r^{2m}}.
\label{eq:poly-optimized}
\end{equation}
The volume fractions at all radii give
\begin{equation}
\lambda_{\Om}^{(m)}
\ge C_{m,d}\,\frac{2m(d+2m)}{d}
\int_0^\infty (1-\Psi_r)\,\frac{\dd r}{r^{2m+1}} .
\label{eq:poly-integral}
\end{equation}
Using the volume fractions at each center gives
\[
\lambda_\Om^{(m)}
\ge C_{m,d}\inf_{x\in\Om}\mathcal L_{2m,0}(\Om,x).
\]
\end{theorem}

\begin{proof}
Fix $r>0$. Apply Theorem~\ref{thm:geometric} with $\alpha=2m$ to
each $x\in\Om$. Since $\psi_r(x)\le\Psi_r$ and $\Phi_{2m,d}$ is
decreasing, we have
\[
M_{\Om}^{(m)}(x)
\ge \frac{\Phi_{2m,d}(\psi_r(x))}{r^{2m}}
\ge \frac{\Phi_{2m,d}(\Psi_r)}{r^{2m}}.
\]
Substituting the last bound into Owen's inequality \eqref{eq:Owen} gives
\[
h_{\Om}^{(m)}[u]
\ge \frac{C_{m,d}}{r^{2m}}\Phi_{2m,d}(\Psi_r)\|u\|_{L^2(\Om)}^2
\qquad (u\in H_0^m(\Om)).
\]
The variational principle gives \eqref{eq:poly-main}.

Evaluating $\Phi_{2m,d}(\Psi_r)$ in the two cases of \eqref{eq:Phi}
gives \eqref{eq:poly-piecewise}. Since $r>0$ was arbitrary, taking
the supremum of the right-hand side of \eqref{eq:poly-main} over all
such $r$ gives \eqref{eq:poly-optimized}. For \eqref{eq:poly-integral},
apply Proposition~\ref{prop:all-radii} with $\alpha=2m$ and $\ell=0$ to
$M_{\Om}^{(m)}$ in \eqref{eq:Owen}. The integral term gives a strictly
larger bound than \eqref{eq:poly-optimized} on bounded centrally
symmetric convex sets other than balls, by
Proposition~\ref{prop:symmetric-convex}.

Since $M_\Om^{(m)}=\mathcal W_{2m,0}$, Theorem~\ref{thm:all-radii-sharp}
gives $M_\Om^{(m)}(x)\ge\mathcal L_{2m,0}(\Om,x)$.
Substituting the infimum of the right-hand side over $x\in\Om$ into
\eqref{eq:Owen} yields the last bound by the same variational principle.
\end{proof}

\begin{corollary}
\label{cor:inradius}
Under the assumptions of Theorem~\ref{thm:polyharmonic}, for every
$q\in(0,1)$,
\begin{equation}
\lambda_{\Om}^{(m)}
\ge \frac{C_{m,d}\Phi_{2m,d}(q)}{(R_{\Om}^{(q)})^{2m}},
\label{eq:inradius-bound}
\end{equation}
where the right-hand side is zero if $R_{\Om}^{(q)}=\infty$.
\end{corollary}

\begin{proof}
Fix $q\in(0,1)$. There is a ball of positive radius contained in
$\Om$. At its center the relative volume at that radius is one,
which gives $R_{\Om}^{(q)}>0$.

Suppose first that $R_{\Om}^{(q)}<\infty$. Each radius for which
$\Psi_r\ge q$ belongs to the set defining $R_{\Om}^{(q)}$. It
follows that $\Psi_r<q$ for every $r>R_{\Om}^{(q)}$. For any such
$r$, the monotonicity of $\Phi_{2m,d}$ and \eqref{eq:poly-main}
give
\[
\lambda_{\Om}^{(m)}\ge \frac{C_{m,d}\Phi_{2m,d}(q)}{r^{2m}}.
\]
Let $r$ decrease to $R_{\Om}^{(q)}$ to obtain
\eqref{eq:inradius-bound}. If $R_{\Om}^{(q)}=\infty$, the
right-hand side of \eqref{eq:inradius-bound} is zero. In this case
the inequality follows from nonnegativity of the quadratic form.
\end{proof}

\subsection{The Dirichlet sub-Laplacian on the Heisenberg group}
\label{subsec:heisenberg-application}

Let $N\ge1$. We use the normalization of
\cite[Subsection~1.3 and Section~6]{FrankLaptevSuragan2026}.
The group law on $\mathbb H^N=\R^{2N}\times\R$ is given by
\[
(z,t)\circ(v,s)=\bigl(z+v,t+s+\tfrac12 z\cdot Jv\bigr),
\qquad
J=\begin{pmatrix}0&I_N\\-I_N&0\end{pmatrix}.
\]
The horizontal vector fields are
$Z_j=\partial_{z_j}-\tfrac12(e_j\cdot Jz)\partial_t$ for
$j=1,\ldots,2N$, where $e_j$ is the $j$th coordinate vector in $\R^{2N}$.
For a nonempty open set $\Om\subset\mathbb H^N$, let $\lambda_\Om^H$
be the bottom of the spectrum of the Dirichlet realization of
$-\sum_{j=1}^{2N}Z_j^2$. The horizontal hyperplane $H(p)$ and the image
$B_r^H(p)$ of the Euclidean ball of radius $r$ under its parametrization
are given by
\[
H(p):=\{p\circ(v,0):v\in\R^{2N}\},
\qquad
B_r^H(p):=\{p\circ(v,0):\abs{v}<r\}.
\]
Write $\sigma$ for the induced $2N$-dimensional Lebesgue measure on
$H(p)$. For $p\in\Om$, set
\[
\widetilde\delta_\omega(p)
:=\inf\{\abs{s}:p\circ(s\omega,0)\notin\Om\},
\qquad \omega\in\Sph^{2N-1},
\]
with $\inf\varnothing=\infty$. The volume fractions in the horizontal
sections are
\begin{equation}
\widetilde\psi_r(p)
:=\frac{\sigma(\Om\cap B_r^H(p))}{\sigma(B_r^H(p))},
\qquad
\widetilde\Psi_r:=\sup_{p\in\Om}\widetilde\psi_r(p).
\label{eq:Heis-density-def}
\end{equation}
The Davies--Hardy distance is
\begin{equation}
\widetilde\delta(p)
:=
\left(
\frac{2N}{|\Sph^{2N-1}|}
\int_{\Sph^{2N-1}}\widetilde\delta_\omega(p)^{-2}\,\dd\omega
\right)^{-1/2},
\label{eq:Heis-Hardy-distance}
\end{equation}
where $\widetilde\delta(p)=\infty$ when the integral is zero.
For geometric Hardy inequalities with other boundary distances on
the Heisenberg group, see Larson~\cite{Larson2016} and Barbatis,
Chatzakou, and Tertikas~\cite{BarbatisChatzakouTertikas2026}.

Fix $p\in\Om$ and set
\[
D_p:=\{v\in\R^{2N}:p\circ(v,0)\in\Om\}.
\]
The set $D_p$ is open and contains the origin. For $p=(z,t)$,
the affine parametrization $v\mapsto p\circ(v,0)$ has derivative
with top block $I_{2N}$ and last row $z^{\mathsf T}J/2$.
Its Gram matrix is $I_{2N}+(J^{\mathsf T}z)(J^{\mathsf T}z)^{\mathsf T}/4$.
Since $J$ is orthogonal, the area Jacobian is
$\sqrt{1+|z|^2/4}$. It is independent of $v$ and cancels in the
volume fraction. Also $s\omega\in D_p$ if and only if
$p\circ(s\omega,0)\in\Om$. Thus
\begin{equation}
\widetilde\psi_r(p)=\frac{|D_p\cap B_r(0)|}{|B_r(0)|},
\qquad
\widetilde\delta_\omega(p)=\delta^{D_p}_\omega(0),
\label{eq:Heis-density}
\end{equation}
where $\delta^{D_p}_\omega$ denotes the directional distance in $D_p$.

\begin{theorem}
\label{thm:heisenberg}
Let $\Om\subset\mathbb H^N$ be nonempty and open. For every $r>0$,
\begin{equation}
\lambda_{\Om}^H
\ge \frac{N}{2r^2}\Phi_{2,2N}(\widetilde\Psi_r).
\end{equation}
In explicit form,
\begin{equation}
\lambda_{\Om}^H
\ge
\begin{cases}
\dfrac{N}{2r^2}\,\widetilde\Psi_r^{-1/N},
& 0<\widetilde\Psi_r\le\dfrac1{N+1},\\[10pt]
\dfrac{(N+1)^{(N+1)/N}}{2r^2}(1-\widetilde\Psi_r),
& \dfrac1{N+1}\le\widetilde\Psi_r\le1.
\end{cases}
\label{eq:Heis-piecewise}
\end{equation}
One also has
\begin{equation}
\lambda_\Om^H\ge(N+1)\int_0^\infty
(1-\widetilde\Psi_r)\,\frac{\dd r}{r^3}.
\label{eq:heis-integral}
\end{equation}
The volume fractions at each center give the further bound
\[
\lambda_\Om^H
\ge\frac N2\inf_{p\in\Om}\mathcal L_{2,0}(D_p,0).
\]
\end{theorem}

\begin{proof}
Fix $r>0$. The form domain of the Dirichlet sub-Laplacian is the
completion of $C_c^\infty(\Om)$ in the norm
\[
\left(\|u\|_{L^2(\Om)}^2+\sum_{j=1}^{2N}\|Z_j u\|_{L^2(\Om)}^2\right)^{1/2}.
\]
For every $u$ in this domain, the Hardy inequality of
Frank--Laptev--Suragan~\cite[Proposition~6.1]{FrankLaptevSuragan2026}
gives
\[
\int_\Om\sum_{j=1}^{2N}\abs{Z_j u(p)}^2\,\dd p
\ge \frac14\int_\Om\abs{u(p)}^2\widetilde\delta(p)^{-2}\,\dd p.
\]
By \eqref{eq:Heis-density}, we may apply Theorem~\ref{thm:geometric}
to $D_p$ at the origin with $d=2N$ and $\alpha=2$. Using the
normalization in \eqref{eq:Heis-Hardy-distance}, the inequality
$\widetilde\psi_r(p)\le\widetilde\Psi_r$, and the monotonicity of
$\Phi_{2,2N}$, we obtain
\begin{align*}
\frac14\widetilde\delta(p)^{-2}
&=\frac{N}{2|\Sph^{2N-1}|}
  \int_{\Sph^{2N-1}}\widetilde\delta_\omega(p)^{-2}\,\dd\omega\\
&\ge \frac{N}{2r^2}\Phi_{2,2N}(\widetilde\psi_r(p))
\ge \frac{N}{2r^2}\Phi_{2,2N}(\widetilde\Psi_r).
\end{align*}
Substituting the last bound into the Hardy inequality and applying
the variational principle proves the first assertion.

To obtain the explicit formula, substitute $\alpha=2$ and $d=2N$
in \eqref{eq:Phi}. The transition point is
$\alpha/(d+\alpha)=1/(N+1)$, and the exponent in the power branch
is $-\alpha/d=-1/N$. The coefficient in the linear branch of
$\Phi_{2,2N}$ is
\[
\frac{d+\alpha}{d}
\left(\frac{d+\alpha}{\alpha}\right)^{\alpha/d}
=\frac{N+1}{N}(N+1)^{1/N}.
\]
Multiplication by $N/2$ gives $(N+1)^{(N+1)/N}/2$.
This proves \eqref{eq:Heis-piecewise}. For \eqref{eq:heis-integral},
apply \eqref{eq:all-radii-pointwise} at the origin of $D_p$ and use
\eqref{eq:Heis-density}. Since $(N/2)k_{2,2N,0}(r)=(N+1)r^{-3}$,
the two volume fractions give the pointwise comparisons
\[
\frac14\widetilde\delta(p)^{-2}
\ge (N+1)\int_0^\infty(1-\widetilde\psi_r(p))\,\frac{\dd r}{r^3}
\ge (N+1)\int_0^\infty(1-\widetilde\Psi_r)\,\frac{\dd r}{r^3}.
\]
For a fixed $p$, the fractions are continuous in $r$ by ball
inclusion and the annulus volume bound. Their supremum
$\widetilde\Psi_r$ is measurable and equals one near zero, since
$D_p$ contains a ball about the origin for any fixed $p\in\Om$.
The integrals are therefore finite. Inserting the last,
center-independent term into the horizontal Hardy inequality and
taking the infimum of the Rayleigh quotient gives
\eqref{eq:heis-integral}.

Alternatively, \eqref{eq:Heis-density} identifies the normalized
spherical average of $\widetilde\delta_\omega(p)^{-2}$ with the
Hardy weight at the origin in $D_p$. Theorem~\ref{thm:all-radii-sharp} bounds
this weight below by $\mathcal L_{2,0}(D_p,0)$.
Multiplying by $N/2$ and taking the infimum over $p\in\Om$ gives
a uniform lower bound for $\widetilde\delta(p)^{-2}/4$. Substitution
in the horizontal Hardy inequality, followed by the Rayleigh
infimum, proves the last assertion.
\end{proof}

\begin{example}
\label{ex:heisenberg-cylinder}
Let $\Om=E\times(-L,L)$, where $L>0$ and $E\subset\R^{2N}$ is a
nonempty bounded centrally symmetric open convex set other than a
ball, with center at the origin. The group law gives, for $p=(z,t)$,
\[
D_p=\left\{v:z+v\in E,\ \left|t+\tfrac12 z\cdot Jv\right|<L\right\}.
\]
In particular, $D_p\subset E-z$ for every $p$ and $D_{(0,0)}=E$.
By \eqref{eq:Heis-density} and the centered comparison in
Proposition~\ref{prop:symmetric-convex},
\[
\widetilde\psi_r(p)\le\frac{|(E-z)\cap B_r(0)|}{|B_r|}
=\frac{|E\cap B_r(z)|}{|B_r|}
\le\frac{|E\cap B_r(0)|}{|B_r|}.
\]
All these inequalities are equalities at $(0,0)$, so
$\widetilde\Psi_r=|E\cap B_r(0)|/|B_r|$ for every $r>0$.
Apply the strict comparison in Proposition~\ref{prop:symmetric-convex}
with $d=2N$, $\alpha=2$, $\ell=0$, and multiply by $N/2$.
Using the zero-shift identities for $k$ and $G$ gives
\[
(N+1)\int_0^\infty (1-\widetilde\Psi_r)\,\frac{\dd r}{r^3}
>\frac N2\sup_{r>0}\frac{\Phi_{2,2N}(\widetilde\Psi_r)}{r^2}.
\]
Thus the integral gives a strictly larger spectral lower bound
on this bounded cylinder.
\end{example}

\section{Comparison with previous bounds}
\label{sec:comparisons}

\subsection{Hardy weights}
\label{subsec:directional-comparison}

For $\ell=0$, the two-branch estimate improves the earlier linear
bound at a fixed radius when the volume fraction is below $q_*$.
After optimization over the radius the two bounds agree. For
$\ell>0$, the sharp bound remains strictly larger than the earlier
shifted estimate after optimization whenever $\Psi_r<1$ at some radius.

With $q=\psi_r(x)$, the zero-shift estimates in
\cite[Lemma~3.1]{FrankLaptevSuragan2026} are
\begin{equation}
\mathcal W_{\alpha,0}(x)\ge\frac{1}{r^\alpha}(q^{-\alpha/d}-1)
\label{eq:old-power}
\end{equation}
and
\begin{equation}
\mathcal W_{\alpha,0}(x)\ge r^{-\alpha}K_{\alpha,d}(1-q).
\label{eq:old-linear}
\end{equation}
For any $\ell\ge0$, the corresponding shifted estimate is
\begin{equation}
\mathcal W_{\alpha,\ell}(x)\ge (r+\ell)^{-\alpha}(1-q).
\label{eq:old-shifted}
\end{equation}

The line $K_{\alpha,d}(1-q)$ is tangent to $q^{-\alpha/d}$ at $q_*$.
At a fixed radius the sharp bound strictly improves both earlier bounds
when $0<q<q_*$. For $q\ge q_*$ it equals the linear bound. On
$[q_*,1]$, the difference between the line and
$q^{-\alpha/d}-1$ is concave, with values $1$ at $q_*$ and $0$ at
$1$. Hence the linear bound also exceeds \eqref{eq:old-power} on
$[q_*,1)$.

The radius change also recovers the pointwise estimate
\eqref{eq:main-geom-intro}. Fix $x$ and set $q=\psi_r(x)$. When
$q\ge q_*$, \eqref{eq:old-linear} at $r$ is the desired bound.
When $q<q_*$, take $s=r(q/q_*)^{1/d}$. Ball inclusion gives
$\psi_s(x)\le(r/s)^d q=q_*$, and \eqref{eq:old-linear} at $s$
bounds the Hardy weight below by
\[
K_{\alpha,d}\frac{1-\psi_s(x)}{s^\alpha}
\ge K_{\alpha,d}\frac{1-q_*}{s^\alpha}
=\frac{q^{-\alpha/d}}{r^\alpha}
=\frac{\Phi_{\alpha,d}(q)}{r^\alpha}.
\]
The convex-envelope argument identifies the equality cases of the
moment problem and thus provides the constructions in
Theorem~\ref{thm:sharpness}.

Proposition~\ref{prop:radius-optimization} gives the same
optimal value for \eqref{eq:old-linear} and the sharp two-branch
bound after taking the supremum over the radius.

To compare the bounds for $\ell>0$, we first strengthen
\eqref{eq:old-shifted}. Use
\eqref{eq:shifted-scalar-linear} with $\ell=0$ and the fact that the
radial volume deficit is at most one. In the notation of
\eqref{eq:K-alpha-d}, this gives the bound in
\cite[Lemma~3.1, Step~2, p.~12]{FrankLaptevSuragan2026},
\[
1-\psi_r(x)\le
\frac{1}{\abs{\Sph^{d-1}}}\int_{\Sph^{d-1}}
\min\left\{1, K_{\alpha,d}^{-1}
\left(\frac{r}{\delta_\omega(x)}\right)^\alpha\right\}\,\dd\omega.
\]
Write $b=K_{\alpha,d}^{-1/\alpha}$. For every $t>0$,
\[
\min\left\{1,\left(\frac{br}{t}\right)^\alpha\right\}
\le \left(\frac{br+\ell}{t+\ell}\right)^\alpha.
\]
For $t\le br$ the right-hand side is at least one. For $t\ge br$,
the inequality $br/t\le(br+\ell)/(t+\ell)$ follows from
$br\,\ell\le t\ell$. Both sides tend to zero as $t\to\infty$, so the
same bound holds in directions of infinite distance. Integrating with
$t=\delta_\omega(x)$ gives
\begin{equation}
\mathcal W_{\alpha,\ell}(x)\ge
\frac{1-\psi_r(x)}{(br+\ell)^\alpha}.
\label{eq:previous-shifted-stronger}
\end{equation}
Since $0<b<1$, this improves \eqref{eq:old-shifted} whenever
$\psi_r(x)<1$.

To compare this coefficient with \eqref{eq:shifted-linear-coef},
use the maximum of $g_r$ at $\ell=0$. From
$C_{\alpha,d,0}(r)=r^{-\alpha}K_{\alpha,d}$ we have
\[
\max_{0\le t\le r} t\left(1-(t/r)^d\right)^{1/\alpha}=br.
\]
For $\ell>0$, use this maximum to bound the term involving $T_r$ and
use $0<T_r<r$ in the remaining term. This yields
\begin{align*}
C_{\alpha,d,\ell}(r)^{-1/\alpha}
&=(T_r+\ell)\left(1-(T_r/r)^d\right)^{1/\alpha}\\
&\le br+\ell\left(1-(T_r/r)^d\right)^{1/\alpha}
<br+\ell.
\end{align*}
Taking inverse powers gives $C_{\alpha,d,\ell}(r)>(br+\ell)^{-\alpha}$
for every positive shift. The two coefficients agree at $\ell=0$.
Thus when $\ell>0$ and $\psi_r(x)<1$, the sharp bound strictly improves
\eqref{eq:old-shifted} and the bound with coefficient
$(br+\ell)^{-\alpha}$.

\begin{proposition}
\label{prop:shifted-strict-comparison}
Under the assumptions of Proposition~\ref{prop:shifted-radius-optimization},
suppose $\ell>0$ and $\Psi_r<1$ for some $r>0$. With
$b=K_{\alpha,d}^{-1/\alpha}$, one has
\begin{equation}
\sup_{r>0}G_{\alpha,d,\ell,r}(\Psi_r)
>
\sup_{r>0}\frac{1-\Psi_r}{(br+\ell)^\alpha}
\ge
\sup_{r>0}\frac{1-\Psi_r}{(r+\ell)^\alpha}.
\label{eq:shifted-strict-comparison}
\end{equation}
If $\Psi_r=1$ for all $r>0$, the three suprema are zero.
\end{proposition}

\begin{proof}
To obtain strict improvement between the suprema, we use a radius
that maximizes the bound with coefficient $(br+\ell)^{-\alpha}$.
Since $\Om$ contains a ball, $\Psi_r=1$ for all sufficiently small
$r$. The function inside the middle supremum in
\eqref{eq:shifted-strict-comparison} therefore vanishes near zero and
tends to zero at infinity since
$1-\Psi_r\le1$. If $\Psi_r<1$ at some radius, the expression is positive
there. Choose a compact interval around that radius so that the value
outside the interval is less than half the positive value. By
\eqref{eq:volume-continuity}, the maximum on the interval is attained at
a radius $r_0$ and is the global maximum. In particular,
$1-\Psi_{r_0}>0$, so the strict
coefficient comparison gives
\[
C_{\alpha,d,\ell}(r_0)(1-\Psi_{r_0})
>\frac{1-\Psi_{r_0}}{(br_0+\ell)^\alpha}.
\]
The supremum of the sharp linear expression is at least its value at
$r_0$. Using \eqref{eq:shifted-radius-optimization} now gives the first
inequality in \eqref{eq:shifted-strict-comparison}. The second follows
from $0<b<1$. If $\Psi_r=1$ at every radius, both linear expressions
and $G_{\alpha,d,\ell,r}(1)$ vanish.
\end{proof}

\subsection{The Robin Laplacian}

Related shifted distance estimates for the Robin $p$-Laplacian on
convex domains are given in \cite{EkholmKovarikLaptev2015}.
Both spectral bounds in \eqref{eq:robin-baselines} follow from
Lemma~\ref{lem:robin-hardy}. For the first we use the geometric
estimate \cite[(3.5)]{FrankLaptevSuragan2026}. For the second we use
\eqref{eq:previous-shifted-stronger}, obtained above by retaining
the coefficient in the proof of that estimate.

\begin{corollary}
\label{prop:robin-comparison}
Under the assumptions of Theorem~\ref{thm:robin}, let
$a=K_{2,d}^{-1/2}$. For every $r>0$, the bottom of the Robin spectrum
has the two lower bounds
\begin{equation}
\frac d4\,\frac{1-\Psi_r}{(r+\ell_\sigma)^2},
\qquad
\frac d4\,\frac{1-\Psi_r}{(ar+\ell_\sigma)^2}.
\label{eq:robin-baselines}
\end{equation}
If $\Psi_r<1$ for some $r>0$, the optimized bound
\eqref{eq:robin-optimized} strictly exceeds the supremum of each
expression in \eqref{eq:robin-baselines}. If $\Psi_r=1$ for every
$r>0$, all three optimized lower bounds are zero.
\end{corollary}

\begin{proof}
Apply \eqref{eq:old-shifted} and its strengthened version
\eqref{eq:previous-shifted-stronger} with $\alpha=2$ and
$\ell=\ell_\sigma$. Each gives a lower bound of the form
$(1-\psi_r(x))/(br+\ell_\sigma)^2$ for the Hardy weight, with
$b=1$ or $b=a$. Since $1-\psi_r(x)\ge1-\Psi_r$, \eqref{eq:robin-hardy}
gives
\[
h_\Om^{\mathrm{R}\sigma}[u]
\ge \frac d4\,\frac{1-\Psi_r}{(br+\ell_\sigma)^2}
\|u\|_{L^2(\Om)}^2,\qquad b\in\{1,a\}.
\]
Taking the infimum of the Rayleigh quotient gives
\eqref{eq:robin-baselines}. Now set $\alpha=2$ and
$\ell=\ell_\sigma>0$ in
Proposition~\ref{prop:shifted-strict-comparison} and multiply
\eqref{eq:shifted-strict-comparison} by $d/4$. Its three expressions
are the optimized bound in \eqref{eq:robin-optimized} and the two
suprema in \eqref{eq:robin-baselines}.
\end{proof}

\begin{example}
\label{ex:robin-disk}
Let $\Om=B_R(0)\subset\R^2$ and take $\sigma=2/(\pi R)$, so that
$\ell_\sigma=R$. At $r=2R$, the volume
fraction is $1/4$. Since $d=\alpha=2$, dividing
\eqref{eq:shifted-tangency} by $r^2$ shows that the ratio
$t_*=T_r/r$ satisfies $4t_*^2+(2R/r)t_*=2$, hence
$4t_*^2+t_*=2$. The left-hand
side at $t=1/2$ is $3/2$, so $1/4<t_*^2$ and the first branch of
\eqref{eq:shifted-profile} applies. Thus
\[
\lambda_{\Om}^{\mathrm{R}\sigma}
\ge \frac12\bigl(R+2R(1/4)^{1/2}\bigr)^{-2}
=\frac{1}{8R^2}.
\]
At the same radius the two bounds in \eqref{eq:robin-baselines}
give $1/(24R^2)$ and $3/(32R^2)$. Indeed,
$K_{2,2}=4$ gives $a=1/2$, so the corresponding denominators
are $(3R)^2$ and $(2R)^2$ and the common numerator is $3/8$.
Here $D=2R$, and the additional term
\eqref{eq:robin-diameter-remainder} is $1/(8R^2)$. Thus the direct
bound with this term is $1/(4R^2)$.
\end{example}

\begin{remark}
On the bounded domains covered by
\cite[Theorem~4.2]{KovarikLaptev2012}, the constant-$\sigma$ Hardy
inequality uses the normalized spherical average with coefficient
$1/4$ and shift $\ell_0=(2\sigma)^{-1}$. Applying
Theorem~\ref{thm:shifted-geometric} to that inequality gives the
spectral lower bound $\frac14\sup_{r>0}G_{2,d,\ell_0,r}(\Psi_r)$.
For $d\ge2$, integrating $|\omega_1|^2<|\omega_1|$ almost everywhere
and using $\int\omega_1^2\,\dd\mu=1/d$ gives $\beta_d>1$. The strict
Cauchy--Schwarz inequality gives $\beta_d<\sqrt d$ because
$|\omega_1|$ varies on the sphere. For each $t\ge0$,
\[
(t+\ell_\sigma)^{-2}
=(t+\beta_d\ell_0)^{-2}
\ge\beta_d^{-2}(t+\ell_0)^{-2}.
\]
The volume constraint in the minimization is the same for both shifts.
Taking the infimum over open sets and centers with the prescribed
volume fraction therefore gives
\[
G_{2,d,\ell_\sigma,r}(q)
\ge\beta_d^{-2}G_{2,d,\ell_0,r}(q).
\]
Including the coefficients of the spherical averages and taking the
supremum over $r>0$ shows that \eqref{eq:robin-optimized} is at least
$d/\beta_d^2>1$ times the optimized bound obtained by applying
the sharp bound to \cite[Theorem~4.2]{KovarikLaptev2012}.
\end{remark}

\begin{remark}
On bounded domains satisfying the strong local Lipschitz condition
in \cite[Assumption~2.1]{JinWeiTangZhan2022}, we can compare with
their Theorem~3.1. At $p=2$, their constants are
$C_2=1/4$ and
\[
B(d,2)=\frac{\sqrt\pi\,\Gamma(d/2+1)}
{\Gamma(3/2)\Gamma(d/2)}=d,
\]
where $\Gamma(3/2)=\sqrt\pi/2$ and
$\Gamma(d/2+1)=(d/2)\Gamma(d/2)$. A constant parameter $\tau>0$ in their
notation gives the coefficient $d\tau$ in front of the boundary
integral and the shift $(2\tau)^{-1}$ in the Hardy weight.
Thus matching the form \eqref{eq:robin-form} requires $\tau=\sigma/d$,
and their inequality becomes
\[
h_{\Om}^{\mathrm{R}\sigma}[u]
\ge \frac d4\int_\Om \mathcal W_{2,d/(2\sigma)}(x)|u(x)|^2\,\dd x .
\]
The same estimate follows from their Theorem~3.2 when the Dirichlet
part of the boundary is empty. The coefficient of the spherical
average in \eqref{eq:robin-hardy} is also $d/4$. Its shift is smaller
because $\beta_d<\sqrt d<d$.

More explicitly, if $D=\operatorname{diam}\Om<\infty$, then
$\delta_\omega(x)\le D/2$, since the nearer endpoint of a line
component is at most half its length away. The function
$t\mapsto (t+d/(2\sigma))/(t+\ell_\sigma)$ has derivative
$(\ell_\sigma-d/(2\sigma))/(t+\ell_\sigma)^2<0$.
It is therefore bounded below on $[0,D/2]$ by its value at $D/2$.
Squaring that lower bound and multiplying by
$(t+d/(2\sigma))^{-2}$ gives an inequality for the two inverse
squares. Averaging it with $t=\delta_\omega(x)$ yields
\[
\mathcal W_{2,\ell_\sigma}(x)
\ge
\left(\frac{D/2+d/(2\sigma)}{D/2+\ell_\sigma}\right)^2
\mathcal W_{2,d/(2\sigma)}(x).
\]
The factor is strictly greater than one on the common
bounded-domain class.
\end{remark}

\begin{example}
Let $\Om=B_R(0)$. Then $\Psi_r=\min\{1,(R/r)^d\}$, since the
intersection is at most as large as the smaller ball and coincident
centers attain this value. At the center $\delta_\omega(0)=R$ in
every direction, so \eqref{eq:shifted-geom} gives
$G_{2,d,\ell_\sigma,r}(\Psi_r)\le(R+\ell_\sigma)^{-2}$. The defining
equation \eqref{eq:shifted-tangency} implies that $T_r\to\infty$ as
$r\to\infty$. For $r>R$ with $T_r\ge R$ we have
$\Psi_r=(R/r)^d\le q_r^*$, and the first branch of
\eqref{eq:shifted-profile} attains this upper bound. Thus the
optimized lower bound in \eqref{eq:robin-optimized} on a ball is
\begin{equation}
\frac{d}{4(R+\ell_\sigma)^2}
=\frac{1}{R^2}\frac{d s^2}{(2s+\beta_d)^2},
\qquad s=\sigma R.
\label{eq:robin-ball}
\end{equation}
For the disk and the parameter in Example~\ref{ex:robin-disk}, this
equals the fixed-radius value $1/(8R^2)$.

Before adding the diameter term, we compare \eqref{eq:robin-ball}
with two bounds in terms of the inradius. On $B_R(0)$,
Kova\v{r}\'{i}k~\cite[Theorem~4.4]{Kovarik2014} gives
$R^{-2}s/[4(1+s)]$, and Proposition~A.2 of Frank, Laptev, and
Suragan~\cite{FrankLaptevSuragan2026} gives
\[
\frac{1}{R^2}\frac{2s^2}{(1+2s)^2}
=\frac{1}{2(R+\ell_0)^2},\qquad \ell_0=(2\sigma)^{-1}.
\]
For bounded convex domains the same inradius bound follows from
Kova\v{r}\'{i}k and Laptev~\cite[Theorem~3.1]{KovarikLaptev2012}.
With $b(x)=\operatorname{dist}(x,\partial\Om)$, their form
inequality has bulk weight
$\tfrac14[(b(x)+\ell_0)^{-2}+(R+\ell_0)^{-2}]$
and an additional nonnegative boundary term. Dropping that term and
using $b(x)\le R$ gives the displayed inradius bound.

For small $s$, Kova\v{r}\'{i}k's bound exceeds \eqref{eq:robin-ball}
because it is of order $s/R^2$, whereas \eqref{eq:robin-ball} is of order
$s^2/R^2$. As $s\to\infty$, the coefficients of $R^{-2}$ tend to
$1/4$ in Kova\v{r}\'{i}k's bound, $1/2$ in Proposition~A.2, and
$d/4$ in \eqref{eq:robin-ball}. Thus for $d\ge3$ our bound exceeds
both inradius bounds for sufficiently large $s$.

On the disk, Proposition~A.2 gives a strictly larger bound than
\eqref{eq:robin-ball} for every $s>0$. The two expressions have
the same numerator $2s^2/R^2$, and $2s+1<2s+\beta_2$ since
$\beta_2=4/\pi>1$. In fact, taking inscribed balls
with radii tending to $R$ gives
$\Psi_r\ge\min\{1,(R/r)^d\}$ on any open set with inradius $R$.
Since $G_{2,d,\ell_\sigma,r}$ is decreasing, the optimized bound
\eqref{eq:robin-optimized} on such a set is at most \eqref{eq:robin-ball}.
By positivity of the kernel, the integral bound \eqref{eq:robin-integral}
is also at most its ball value, which equals \eqref{eq:robin-ball} by
Proposition~\ref{prop:symmetric-convex}. Thus, on planar convex domains
of finite inradius, the inradius bound $1/[2(R+\ell_0)^2]$ strictly
exceeds the optimized single-radius bound and the integral bound before
the diameter term is added.

The diameter term \eqref{eq:robin-diameter-remainder} changes the
comparison on balls. Here $D=2R$, so it doubles
\eqref{eq:robin-ball} and gives
\[
\lambda_{B_R}^{\mathrm{R}\sigma}\ge
\frac{d}{2(R+\ell_\sigma)^2}.
\]
The ratio of this bound to $1/[2(R+\ell_0)^2]$ is
$d(R+\ell_0)^2/(R+\beta_d\ell_0)^2\ge d/\beta_d^2>1$.
Thus on balls in every dimension $d\ge2$, adding the diameter term
gives a strictly larger bound than the inradius bound of Proposition~A.2.

The Robin Faber--Krahn bound of \cite[Theorem~1.1]{Daners2006} uses
the first eigenvalue of the ball of the same volume. On $B_R(0)$ it
gives the exact eigenvalue, which exceeds even the bound with the
diameter term.
Indeed, solving $|x+s\omega|=R$ for the two endpoints gives
\[
\delta_\omega(x)=\sqrt{(x\cdot\omega)^2+R^2-|x|^2}
-|x\cdot\omega|\le\sqrt{R^2-|x|^2}.
\]
For $x\ne0$ the last expression is smaller than $R$, so the
Hardy weight is strictly greater than $(R+\ell_\sigma)^{-2}$
except at the center. A nonzero first eigenfunction has positive
squared $L^2$ mass away from the center. Applying the bounded-domain
form estimate in Lemma~\ref{lem:robin-hardy} to that eigenfunction
therefore gives the strict inequality.
\end{example}

\subsection{The Dirichlet Laplacian and polyharmonic operators}

\begin{remark}
Inserting \eqref{eq:old-power} and \eqref{eq:old-linear} into
\eqref{eq:Owen} gives the power and linear spectral
bounds of \cite[Theorem~2.3]{FrankLaptevSuragan2026}. The fixed-radius
comparison in Subsection~\ref{subsec:directional-comparison} with
$\alpha=2m$ shows that \eqref{eq:poly-piecewise} strictly improves
both bounds for $0<\Psi_r<2m/(d+2m)$ and agrees with the linear bound
for $\Psi_r\ge 2m/(d+2m)$.

By Proposition~\ref{prop:radius-optimization} and the identity
$K_{2m,d}=c_{m,d}$, the right-hand side of \eqref{eq:poly-optimized}
equals the optimized linear estimate
\begin{equation}
C_{m,d}c_{m,d}\sup_{r>0}\frac{1-\Psi_r}{r^{2m}}.
\label{eq:poly-linear-optimized}
\end{equation}

For the Dirichlet Laplacian ($m=1$), the earlier estimate of
Frank and Larson~\cite[Theorem~2]{FrankLarson2021}, which also follows
from Davies~\cite[Theorem~1.5.11]{Davies1989}, has right-hand side
$d(1-\Psi_r)/(4r^2)$. By \eqref{eq:Phi}, the right-hand side of
\eqref{eq:poly-main} is at least $K_{2,d}$ times this estimate at
any fixed radius. Equation~\eqref{eq:poly-linear-optimized} gives
exactly the factor $K_{2,d}$ after optimizing both estimates over $r$.

\end{remark}

\begin{example}
For the planar biharmonic operator ($d=m=2$), the constants and the
threshold are
\[
C_{2,2}=\frac32,\qquad
\frac{2m}{d+2m}=\frac23,\qquad
c_{2,2}=\frac{27}{4}.
\]
Consequently, the bound is
\begin{equation}
\lambda_{\Om}^{(2)}
\ge \frac{1}{r^4}
\begin{cases}
\dfrac32\,\Psi_r^{-2}, & 0<\Psi_r\le\dfrac23,\\[8pt]
\dfrac{81}{8}(1-\Psi_r), & \dfrac23\le\Psi_r\le1.
\end{cases}
\end{equation}
At $\Psi_r=1/2$ this gives
\[
\lambda_{\Om}^{(2)}\ge\frac{6}{r^4}.
\]
At the same radius, the spectral bounds from \eqref{eq:old-power}
and \eqref{eq:old-linear} give $9/(2r^4)$ and $81/(16r^4)$,
respectively. Both are strictly smaller than $6/r^4$.
Set $s=(\sqrt3/2)r$. Then \eqref{eq:density-scale} gives
$\Psi_s\le2/3$, so the earlier linear estimate at $s$ gives
\[
\lambda_{\Om}^{(2)}
\ge \frac{81}{8s^4}(1-\Psi_s)
\ge \frac{27}{8s^4}
= \frac6{r^4}.
\]
\end{example}

For the comparison with Lieb's lower bound~\cite{Lieb1983}, let
$\Lambda_B=\lambda_B^{\mathrm D}$ be the first Dirichlet eigenvalue of
the unit ball $B$. Lieb takes the supremum of the ball fractions
over all centers in $\R^d$, whereas $\Psi_r$ in \eqref{eq:poly-main}
uses centers in $\Om$. We first compare the two fixed-radius
expressions at the same $r$ and the same fraction $q\in(0,1)$.
The Lieb expression $\Lambda_B r^{-2}(q^{-2/d}-1)$ is then strictly
larger than $d\Phi_{2,d}(q)/(4r^2)$ in \eqref{eq:poly-main} for
$m=1$. To see this, put $p=2/d$ and
$q_*=p/(1+p)$. For $q\ge q_*$, convexity gives
$q^{-p}-1\ge p(1-q)$, so
$\Phi_{2,d}(q)/(q^{-p}-1)\le K_{2,d}/p$. For $q\le q_*$, the ratio
is $(1-q^p)^{-1}$. It increases up to its value at $q_*$, where the
two branches agree, so the same bound holds. By concavity of
$t\mapsto t^p$,
\[
\frac{dK_{2,d}}{4p}
=\frac{d^2}{8}(1+p)(1+1/p)^p\le\frac{d(d+2)}{4}.
\]
On $B$, the positive function $v(x)=(1-|x|^2)e^{-d|x|^2/4}$ satisfies
\[
\frac{-\Delta v(x)}{v(x)}
=2d+\frac{d^2}{2}-\frac{d^2}{4}|x|^2\ge 2d+\frac{d^2}{4}.
\]
For $u\in C_c^\infty(B)$, expansion and integration by parts give
\[
\int_B|\nabla u|^2\,\dd x
=\int_B v^2|\nabla(u/v)|^2\,\dd x
+\int_B\frac{-\Delta v}{v}|u|^2\,\dd x.
\]
The variational principle yields $\Lambda_B\ge 2d+d^2/4$,
which exceeds $d(d+2)/4$. This proves the strict comparison at a
common $q$.

We now compare the optimized Lieb bound with the integral bound
\eqref{eq:poly-integral}. On the balls and boxes considered below,
the two choices of centers give the same maximal volume fraction
at every radius.

\begin{proposition}
\label{prop:lieb-optimized}
For a nonempty open set $\Om\subset\R^d$, $d\ge2$, write
\begin{align*}
\mathcal B_{\rm Lieb}(\Om)
&=\Lambda_B\sup_{r>0}\frac{\widehat\Psi_r^{-2/d}-1}{r^2},
&
\widehat\Psi_r&=\sup_{x\in\R^d}\frac{|\Om\cap B_r(x)|}{|B_r|},\\
\mathcal B_{\rm int}(\Om)
&=\frac{d+2}{2}\int_0^\infty (1-\Psi_r)\,\frac{\dd r}{r^3}.
\end{align*}
These are the optimized Lieb bound and the right-hand side of
\eqref{eq:poly-integral} for $m=1$, respectively.

For a ball $B_R$ of radius $R>0$, one has
\[
\mathcal B_{\rm Lieb}(B_R)=\Lambda_B R^{-2}>
\frac{d}{4R^2}=\mathcal B_{\rm int}(B_R).
\]
For all sufficiently large $d$ and, for each such $d$, all sufficiently
large $L$, the box $\Om_L=(-1,1)\times(-L,L)^{d-1}$ satisfies
\begin{equation}
\mathcal B_{\rm Lieb}(\Om_L)<\frac14\le\mathcal B_{\rm int}(\Om_L).
\label{eq:lieb-box}
\end{equation}
\end{proposition}

\begin{proof}
For a bounded connected open set $A$ and $r>0$, combining Lieb's
intersection theorem~\cite{Lieb1983} with the Faber--Krahn inequality
gives a center $x$ such that
\[
\Lambda_B r^{-2}\widehat\Psi_r(A)^{-2/d}
\le \lambda_{A\cap B_r(x)}^{\mathrm D}
< \lambda_A^{\mathrm D}+\Lambda_B r^{-2}.
\]
For a bounded open set with several components, apply the resulting
lower bound to each component $C$. We have
$\widehat\Psi_r(C)\le\widehat\Psi_r(A)$ and
$\lambda_A^{\mathrm D}=\inf_C\lambda_C^{\mathrm D}$ by splitting the
Dirichlet energy and squared norm over the components. For general
$\Om$, take $\Om_R=\Om\cap B_R(0)$. Then
$\widehat\Psi_r(\Om_R)\le\widehat\Psi_r(\Om)$ and
$\lambda_{\Om_R}^{\mathrm D}\downarrow\lambda_\Om^{\mathrm D}$, since
every compactly supported smooth test function in $\Om$ is supported
in $\Om_R$ for all sufficiently large $R$. Passing to the limit and
taking the supremum over $r>0$ give
$\lambda_\Om^{\mathrm D}\ge\mathcal B_{\rm Lieb}(\Om)$.

For the ball, $\widehat\Psi_r=\Psi_r=\min\{1,(R/r)^d\}$.
The expression in the Lieb supremum is zero for $r\le R$ and equals
$R^{-2}-r^{-2}$ for $r>R$. The integral equals
\[
\frac{d+2}{2}\int_R^\infty \left(1-(R/r)^d\right)\frac{\dd r}{r^3}
=\frac{d}{4R^2}.
\]
The strict inequality follows from the lower bound for $\Lambda_B$
proved above.

The integral bound on the slab $S_d=(-1,1)\times\R^{d-1}$ is $1/4$.
We show that Lieb's bound is smaller in high dimensions, first on
the slab and then on sufficiently long boxes. The centered-volume argument in
Proposition~\ref{prop:symmetric-convex} applies to both $S_d$ and each
box. Even for the unbounded slab the sets in the Brunn--Minkowski
inequality are its intersections with balls, hence bounded.
Thus the origin maximizes the intersection volume over all centers,
and $\widehat\Psi_r=\Psi_r=\psi_r(0)$ for every $r$.
At the origin of the slab, $\delta_\omega(0)=1/|\omega_1|$, with
infinite distance if $\omega_1=0$. A ray is outside the slab
precisely at radii $r\ge\delta_\omega(0)$ when the distance is
finite. Hence equality holds in \eqref{eq:spherical-deficit} at the
origin and the integral formula in \eqref{eq:all-radii-pointwise}
gives
\[
\mathcal B_{\rm int}(S_d)=\frac d4\mathcal W_{2,0}(0)
=\frac d4\int_{\Sph^{d-1}}\omega_1^2\,\dd\mu(\omega)=\frac14.
\]

Let $X$ be uniformly distributed in the unit ball of $\R^d$ and set
$q_d(s)=\mathbb P(|X_1|\le s)$. For $r\ge1$ the fraction at the
origin of $S_d$ is $q_d(1/r)$. For $r<1$ it is one.
For $0<s<1$, let $g_s(t)$ be the conditional probability of
$|X_2|\le s$ given $|X_1|=t$. The conditional law of the remaining
coordinates is uniform on a $(d-1)$-dimensional ball of radius
$\sqrt{1-t^2}$. Scaling this ball to the unit ball shows that
$g_s$ is nondecreasing. Its average on $\{|X_1|\le s\}$ is
therefore at most its average over all $|X_1|$. By symmetry
the latter average is $\mathbb P(|X_2|\le s)=q_d(s)$.
Consequently,
\[
\mathbb P(|X_1|\le s,\ |X_2|\le s)
=q_d(s)\mathbb E[g_s(|X_1|)\mid |X_1|\le s]\le q_d(s)^2.
\]
The disk $X_1^2+X_2^2\le s^2$ lies inside the square in the event
on the left. The two-coordinate density is proportional to
$(1-x_1^2-x_2^2)^{(d-2)/2}$ on the disk, by taking sections of the
unit ball. Thus the disk of radius $s$ has probability
\[
\frac{\int_0^s \rho(1-\rho^2)^{(d-2)/2}\,\dd \rho}
{\int_0^1 \rho(1-\rho^2)^{(d-2)/2}\,\dd \rho}
=1-(1-s^2)^{d/2}.
\]
Using $1-s^2\le e^{-s^2}$ gives, also at $s=1$,
\[
q_d(s)^2\ge 1-(1-s^2)^{d/2}\ge 1-e^{-d s^2/2}.
\]
For $r\ge1$, put $s=1/r$ and $u=d s^2/2$. The expression under the
supremum defining $\mathcal B_{\rm Lieb}(S_d)$ satisfies
\begin{align*}
\Lambda_B r^{-2}(\widehat\Psi_r^{-2/d}-1)
&=\Lambda_Bs^2(q_d(s)^{-2/d}-1)\\
&\le\frac{2\Lambda_B}{d^2}f_d(u),
\qquad f_d(u)=du\left[(1-e^{-u})^{-1/d}-1\right].
\end{align*}
For $r<1$ the expression is zero. Extending the supremum on the
right to all $u>0$ gives
\begin{equation}
\mathcal B_{\rm Lieb}(S_d)
\le\frac{2\Lambda_B}{d^2}\sup_{u>0} f_d(u).
\label{eq:lieb-slab}
\end{equation}

We next take $d\to\infty$. With $a(u)=-\log(1-e^{-u})$ we have
$f_d(u)=du(e^{a(u)/d}-1)$, which converges to $u a(u)$ uniformly
on every compact subinterval of $(0,\infty)$. Moreover,
$0\le f_d(u)\le ua(u)e^{a(u)/2}$ for $d\ge2$, since
$e^z-1\le ze^z$ for $z\ge0$. The upper bound tends to zero at
both ends of $(0,\infty)$. As $u\downarrow0$ it is
$O(\sqrt u\,(1+|\log u|))$, and as $u\to\infty$ it is
$O(ue^{-u})$. Given $\varepsilon>0$, choose a compact subinterval
outside which the upper bound is less than $\varepsilon$.
The limiting function $ua(u)$ obeys the same bound there.
On the compact subinterval $f_d\to ua(u)$ uniformly, so the
difference of the suprema tends to zero. Substituting $b=e^{-u}$ gives
$u a(u)=\log b\log(1-b)$. Its derivative as a function of $b$ has
the sign of $(1-b)\log(1-b)-b\log b$. This difference vanishes at
$0$ and $1/2$ and has negative second derivative on $(0,1/2)$.
It is positive on that interval and is antisymmetric about $1/2$.
Thus $\sup_{u>0}u a(u)=(\log 2)^2$.

The bound $\Lambda_B\ge2d+d^2/4$ gives
$\liminf_{d\to\infty}\Lambda_B/d^2\ge1/4$. For the upper bound, let
$d\ge3$, fix $0<h<1$, and take
$v(\rho)=\sin(\pi(\rho-1+h)/h)$ on $(1-h,1)$ and zero on
$(0,1-h)$. The radial trial function $|x|^{-(d-1)/2}v(|x|)$
is zero near the origin and has zero trace on the unit sphere,
so it belongs to $H_0^1(B)$. Its squared norm and energy, after
canceling the sphere area, are $\int_0^1|v|^2\,\dd\rho$ and
$\int_0^1|v'-cv/\rho|^2\,\dd\rho$, where $c=(d-1)/2$.
In the latter integral the cross term is
$-c\int_0^1 (|v|^2)'/\rho\,\dd\rho
=-c\int_0^1 |v|^2/\rho^2\,\dd\rho$,
since $v$ vanishes at both endpoints of the layer. The sine has
$\int_0^1|v'|^2\,\dd\rho=(\pi/h)^2\int_0^1|v|^2\,\dd\rho$.
Thus the quotient is
\[
\frac{\int_0^1
\left(|v'|^2+\frac{(d-1)(d-3)}{4\rho^2}|v|^2\right)\,\dd\rho}
     {\int_0^1 |v|^2\,\dd\rho}
\le \frac{\pi^2}{h^2}+\frac{(d-1)(d-3)}{4(1-h)^2}.
\]
Choosing $h=d^{-1/2}$ and dividing by $d^2$ gives
$\limsup_{d\to\infty}\Lambda_B/d^2\le1/4$. Thus
$\Lambda_B/d^2\to1/4$, and \eqref{eq:lieb-slab} yields
\[
\limsup_{d\to\infty}\mathcal B_{\rm Lieb}(S_d)
\le\frac{(\log 2)^2}{2}<\frac14.
\]
For the strict inequality, the first four terms of the
exponential series at $7/10$ sum to more than $2$.
Thus $\log 2<7/10$ and $(\log 2)^2/2<49/200<1/4$.

Fix a dimension in which $\mathcal B_{\rm Lieb}(S_d)<1/4$.
For $r\le L$, every transverse coordinate of a point of $B_r(0)$
has absolute value less than $L$, so the intersections of $\Om_L$
and $S_d$ with $B_r(0)$ coincide. Thus the Lieb
expressions agree at those radii. For $r>L$, the intersection
$\Om_L\cap B_r(0)$ contains $\Om_L\cap B_L(0)$. For $L\ge2$ this
intersection contains
$(-1,1)\times(-L/(2\sqrt{d-1}),L/(2\sqrt{d-1}))^{d-1}$,
since points therein have squared norm less than $1+L^2/4<L^2$.
Its volume is therefore at least $2(L/\sqrt{d-1})^{d-1}$.
Consequently, for $r>L$,
\[
\frac{\Lambda_B}{r^2}(\widehat\Psi_r^{-2/d}-1)
\le \Lambda_B\left(\frac{|B|}{|\Om_L\cap B_L(0)|}\right)^{2/d}
=O_d(L^{-2(d-1)/d}).
\]
The constant is independent of $r>L$. For large $L$ this bound
is below $1/4$, while the supremum for $r\le L$ is bounded above
by $\mathcal B_{\rm Lieb}(S_d)<1/4$. Hence the supremum over
all radii satisfies $\mathcal B_{\rm Lieb}(\Om_L)<1/4$.
Since $\Om_L\subset S_d$, its volume fractions are at most those
of the slab and $\mathcal B_{\rm int}(\Om_L)\ge1/4$. This proves
\eqref{eq:lieb-box}.
\end{proof}

\subsection{The Dirichlet sub-Laplacian on the Heisenberg group}

We use the volume fractions in the horizontal sections from
Subsection~\ref{subsec:heisenberg-application}.
Applying \eqref{eq:old-power} and \eqref{eq:old-linear} to
the horizontal sections with $d=2N$ and $\alpha=2$ in the proof of
Theorem~\ref{thm:heisenberg} gives the power and linear spectral
bounds of \cite[Theorem~2.5]{FrankLaptevSuragan2026}.
The fixed-radius comparison in Subsection~\ref{subsec:directional-comparison}
shows that \eqref{eq:Heis-piecewise} strictly improves both bounds for
$\widetilde\Psi_r<1/(N+1)$ and agrees with the linear bound for
$\widetilde\Psi_r\ge1/(N+1)$. To compare the optimized bounds, note that
$B_s^H(p)\subset B_r^H(p)$ for $0<s\le r$ and that
$\sigma(B_r^H(p))$ is proportional to $r^{2N}$. Hence
\[
s^{2N}\widetilde\psi_s(p)\le r^{2N}\widetilde\psi_r(p).
\]
Taking the supremum over $p\in\Om$ gives the same scaling inequality
for $\widetilde\Psi_r$. The argument of
Proposition~\ref{prop:radius-optimization} with $d=2N$ and $\alpha=2$
then identifies the optimized value with that of the linear bound.

On Heisenberg domains covered by
\cite[Proposition~2, (6)]{PrandiRizziSeri2019}, the Hardy weight is
$N/2$ times the normalized spherical average of the inverse squares
of the one-sided exit lengths
\[
\widetilde\delta_\omega^+(p)
:=\inf\{s>0:p\circ(s\omega,0)\notin\Om\}.
\]
Here $\inf\varnothing=\infty$. For $0\le s<\widetilde\delta_\omega^+(p)$
one has $s\omega\in D_p$.
Integrating these ray segments in polar coordinates and using
\eqref{eq:Heis-density} gives
\[
\widetilde\psi_r(p)\ge
\frac{1}{|\Sph^{2N-1}|}\int_{\Sph^{2N-1}}
\min\left\{1,\left(\frac{\widetilde\delta_\omega^+(p)}{r}\right)^{2N}\right\}
\,\dd\omega.
\]
Apply Lemma~\ref{lem:convex-envelope} with $d=2N$, $\alpha=2$,
and $Y(\omega)=\widetilde\delta_\omega^+(p)/r$. The preceding
inequality gives the moment constraint with $q=\widetilde\psi_r(p)$.
Dividing the conclusion of the lemma by $r^2$ yields
\[
\frac{1}{|\Sph^{2N-1}|}\int_{\Sph^{2N-1}}
   \bigl(\widetilde\delta_\omega^+(p)\bigr)^{-2}\,\dd\omega
\ge r^{-2}\Phi_{2,2N}\bigl(\widetilde\psi_r(p)\bigr).
\]
The variational argument in the proof of Theorem~\ref{thm:heisenberg}
now gives \eqref{eq:Heis-piecewise} on this class of domains.

The spectral estimate in \cite[Corollary~6]{PrandiRizziSeri2019}
uses the horizontal diameter, the supremum of the lengths of straight
horizontal segments contained in the domain. The bounds above use
volume fractions in the horizontal sections.

\end{document}